\documentclass[11pt]{amsart}
\usepackage[a4paper,margin=2.6cm]{geometry}
\usepackage{amsmath,amssymb,amsthm,mathtools,booktabs,enumitem}
\usepackage{xcolor,graphicx}
\usepackage[colorlinks=true,linkcolor=blue!50!black,citecolor=blue!50!black]{hyperref}

\newtheorem{theorem}{Theorem}[section]
\newtheorem{proposition}[theorem]{Proposition}
\newtheorem{lemma}[theorem]{Lemma}
\newtheorem{corollary}[theorem]{Corollary}
\theoremstyle{definition}
\newtheorem{definition}[theorem]{Definition}
\newtheorem{remark}[theorem]{Remark}
\numberwithin{equation}{section}

\newcommand{\M}{\mathbb{M}^n(\varepsilon)}
\newcommand{\E}{\mathbb{E}^{n+1}_\varepsilon}
\newcommand{\s}{\mathrm{s}}
\newcommand{\cc}{\mathrm{c}}
\newcommand{\la}{\langle}
\newcommand{\ra}{\rangle}
\newcommand{\R}{\mathbb{R}}
\newcommand{\Hy}{\mathbb{H}}
\newcommand{\Sp}{\mathbb{S}}
\newcommand{\U}{\mathbb{U}}

\title[Invariant solitons of the IMCF in $\Sp^n$ and $\Hy^n$]{Invariant solitons of the inverse mean curvature flow\\ in the sphere and in hyperbolic space}
\author[H. F. S. dos Reis]{Hiuri Fellipe Santos dos Reis}
\address{Instituto de Matem\'atica e Estat\'istica, Universidade Federal de Goi\'as, 74690-900, Goi\^ania--GO, Brazil}
\email{hiuri\_reis@ufg.br}
\author[L. R. Adriano]{Levi Rosa Adriano}
\address{Instituto de Matem\'atica e Estat\'istica, Universidade Federal de Goi\'as, 74690-900, Goi\^ania--GO, Brazil}
\author[M. C. F. Furtado]{Miriam Cristina Ferreira Furtado}
\address{Instituto de Matem\'atica e Estat\'istica, Universidade Federal de Goi\'as, 74690-900, Goi\^ania--GO, Brazil}


\begin{document}

\begin{abstract}
We study hypersurfaces of the sphere $\Sp^n$ and of hyperbolic space $\Hy^n$, $n\ge3$, that move under the inverse mean curvature flow by a one-parameter group of isometries. No closed hypersurface has this property, and a soliton with constant mean curvature is an open subset of a horosphere. For the rotation hypersurfaces of do Carmo and Dajczer, of spherical, hyperbolic and parabolic type, the soliton equation reduces to a single ordinary differential equation, which we analyse in every dimension. In $\Hy^n$ the horosphere appears at the critical speed $\lambda=1$ and separates entire bowls from solitons that end at a singular $(n-2)$-sphere; each class has exactly one horospherical end, and it is the threshold between ends asymptotic to totally geodesic hypersurfaces and ends at a singular orbit. At an axis of codimension two we obtain a complete trichotomy, which changes at the speed $\lambda=2(n-1)$ from logarithmic spirals to cuspidal edges. In $\Sp^n$ no rotational soliton through the axis is complete, and the radius of its singular sphere satisfies $\frac\pi2-r_+=\frac{K_n}\lambda\big(1-2(n-1)^2\lambda^{-2}\log\lambda\big)+O(\lambda^{-3})$, with $K_n=\sqrt\pi\,\Gamma(\frac{n+1}2)/\Gamma(\frac n2)$. For $n=3$, the commuting pairs of Killing fields give two further classes, helicoidal and Hopf-invariant surfaces; the latter are the rotators of the inverse curve shortening flow on the round two-sphere.
\end{abstract}

\maketitle

\section{Introduction}

The inverse mean curvature flow (IMCF) moves a hypersurface in the direction of its unit normal with speed equal to the inverse of its mean curvature. In contrast with the mean curvature flow, it is an expanding flow. Its weak formulation by Huisken and Ilmanen \cite{HI} led to the proof of the Riemannian Penrose inequality. In hyperbolic space, Gerhardt \cite{Ge} proved that star-shaped mean-convex hypersurfaces become umbilic at infinity under the IMCF. Hung and Wang \cite{HW} showed that the rescaled flow does not, in general, converge to a round sphere. This makes the IMCF in non-Euclidean space forms a natural place to look for explicit solutions.

Self-similar solutions play a central role in the study of curvature flows, since they model singularities and the long-time behaviour. In the Euclidean setting, solitons of the IMCF were studied by Drugan, Lee and Wheeler \cite{DLW}, by Castro and Lerma \cite{CL} and by Kim and Pyo \cite{KP}, among others; see also \cite{Fu} for graphical translating solitons in a product setting. In a Riemannian manifold, the natural analogue of a self-similar solution is a soliton with respect to a Killing field: a solution that evolves by a one-parameter group of isometries of the ambient space. For the mean curvature flow such solitons were introduced by Hungerb\"uhler and Smoczyk \cite{HS}, and recently studied in the sphere, for the Hopf field, by Magliaro, Mari, Roing and Savas-Halilaj \cite{MMRS}. In hyperbolic three-space, de Lima, Ramos and dos Santos \cite{LRS} studied translators and rotators of the mean curvature flow. They constructed translating catenoids, grim reapers and complete helicoidal rotators, proved a tangency principle, and showed that horospheres are the only solitons with nonzero constant mean curvature. Rotational translators in $\Hy^{n+1}$ were then constructed by de Lima and Ramos \cite{dLR}, and translators with respect to a parabolic Killing field of $\Hy^3$ were classified by Bueno and L\'opez \cite{BL}. Our rotational translators and helicoidal solitons are the IMCF counterparts of their examples, and Proposition~\ref{prop:cmc} below is the IMCF version, in every dimension, of their characterisation of horospheres.

For the inverse flow in non-flat space forms, explicit solutions are scarce. Alencar and Tenenblat \cite{AT} described the solutions that move by parallel hypersurfaces, which are the isoparametric ones. In dimension two, Krznari\'c and L\'opez \cite{KL} recently studied the solitons of the inverse curve shortening flow in the hyperbolic plane, where horocycles also play a distinguished role. The present paper deals with hypersurfaces of dimension at least two, in both $\Sp^n$ and $\Hy^n$; here the orbits of the symmetry group have positive dimension, and the dimension enters the critical speeds and the asymptotic constants.

In this paper we consider hypersurfaces $\Sigma^{n-1}\subset\M$, where $\M$ is $\Sp^n$ ($\varepsilon=1$) or $\Hy^n$ ($\varepsilon=-1$) and $n\ge3$, that are IMCF solitons with respect to a Killing field $X$. We first show that there is no closed IMCF soliton (Proposition~\ref{prop:noclosed}), and that a soliton with constant mean curvature, for an arbitrary Killing field, is an open subset of a horosphere (Proposition~\ref{prop:cmc}). The proof of the latter uses that the IMCF of such a soliton moves by parallel hypersurfaces.

We then consider rotation hypersurfaces in the sense of do Carmo and Dajczer \cite{dCD83}: hypersurfaces invariant under the group $G$ of isometries that fix a $2$-plane of $\E$ pointwise. In $\Hy^n$, according to the causal character of the plane, there are three types, spherical, hyperbolic and parabolic; in $\Sp^n$ there is only one, which we also call spherical. In each case $G$ has orbits of codimension two, and there is a Killing field $X_*$ commuting with $G$. We prove that, for all of them and in every dimension, the soliton equation with respect to $X=\lambda X_*$ reduces to a single ordinary differential equation for the angle $\sigma$ between the profile curve and the radial direction (Theorem~\ref{thm:master}). The mean curvature is the sum of the geodesic curvature of the profile and of the logarithmic derivative of the volume of the orbits. We then describe the solutions for each type. Theorems~\ref{thm:transl}--\ref{thm:parab} give the complete qualitative picture: every end of every solution is classified.

The dimension enters the results in a transparent way:
\begin{itemize}[nosep]
\item the critical speed at which the horosphere appears among the rotational translators of $\Hy^n$ is $\lambda=1$ for every $n$ (Theorem~\ref{thm:transl});
\item every class in $\Hy^n$ has exactly one \emph{horospherical end}, with $H\to1$. It is a threshold: the outer ends on one side of it are complete and asymptotic to totally geodesic hypersurfaces at the rate $O(e^{-2r})$, with $\sin\sigma=\frac{2(n-1)}{(n-2)\lambda}e^{-r}(1+o(1))$, and those on the other side stop at a singular orbit (Theorem~\ref{thm:outer});
\item at an axis of codimension two (the totally geodesic $\Hy^{n-2}$ of the hyperbolic type, or the far axis $\Sp^{n-2}$ in $\Sp^n$) the behaviour changes at $\lambda=2(n-1)$. For the hyperbolic type we obtain a complete trichotomy for every $\lambda$: above the threshold speed the solitons reach the axis along logarithmic spirals of two possible angles, one of them being realised by exactly one solution; below it they reach the axis only as fronts, through infinitely many cuspidal edges (Theorems~\ref{thm:hyp} and~\ref{thm:sphere});
\item the radius of the singular sphere of the rotational caps satisfies $r_+=c_0(n)\lambda+O(\lambda^3)$ as $\lambda\to0$, and, in $\Sp^n$,
\[
\frac\pi2-r_+=\frac{K_n}\lambda-\frac{2(n-1)^2K_n\log\lambda}{\lambda^3}+O(\lambda^{-3})\qquad(\lambda\to\infty),
\]
with the Wallis-type constant $K_n=(n-1)\int_0^{\pi/2}\cos^{n-2}t\,dt$ (Theorem~\ref{thm:rplus}). The proof rests on a flux identity (Lemma~\ref{lem:flux}), which says that the volume of the cap is proportional to the flux of the Killing field through its singular sphere, and on a second, Clairaut-type identity that produces the logarithmic term.
\end{itemize}
Every singular orbit is a cuspidal edge of a wave front (Theorem~\ref{thm:fronts}), so the solitons continue naturally beyond the singular set.

For $n=3$ there is more symmetry. A rotation surface is invariant under a one-parameter group $Y$, and the classical classification of commuting pairs $(X,Y)$ (Proposition~\ref{prop:pairs}) adds two classes: helicoidal surfaces and, in $\Sp^3$, Hopf-invariant surfaces. For all of them the soliton equation reduces to one equation whose drift term does not depend on the pitch of $Y$ (Theorem~\ref{thm:master3}). We obtain the three-dimensional results as corollaries of the general theorems (Corollaries~\ref{cor:transl3}--\ref{cor:parab3}), and we complete them with the helicoidal and Hopf classes (Theorems~\ref{thm:helic} and~\ref{thm:hopf}). Through the Hopf fibration, the Hopf-invariant solitons are the rotators of the inverse curve shortening flow on $\Sp^2(\frac12)$. For $n\ge4$ there is no such class in $\Sp^n$ (Remark~\ref{rem:rank}).

The paper is organised as follows. Our results are stated in Section~\ref{sec:main}. Section~\ref{sec:prelim} contains the general properties of solitons, Section~\ref{sec:master} the reduction to an ordinary differential equation and the structure of the singular orbits, Section~\ref{sec:rot} the analysis of the spherical and hyperbolic types, Section~\ref{sec:parab} the parabolic type, and Section~\ref{sec:three} the three-dimensional classes. All statements are proved analytically. Section~\ref{sec:numerics} contains numerical illustrations only.

\section{Main results}\label{sec:main}

\subsection{Setting}
Let $n\ge3$. We denote by $\M$ the space form of dimension $n$ and constant sectional curvature $\varepsilon\in\{1,-1\}$, realised as the quadric
\[
\M=\{x\in\E:\ \la x,x\ra_\varepsilon=\varepsilon\}\quad (x_0>0 \text{ if } \varepsilon=-1),
\]
where $\E=(\R^{n+1},\la\cdot,\cdot\ra_\varepsilon)$ and $\la x,y\ra_\varepsilon=\varepsilon x_0y_0+\sum_{i=1}^nx_iy_i$. As in \cite{RT}, we use the notation
\begin{equation}\label{eq:cs}
\cc(\xi)=\begin{cases}\cos\xi,&\varepsilon=1,\\ \cosh\xi,&\varepsilon=-1,\end{cases}\qquad
\s(\xi)=\begin{cases}\sin\xi,&\varepsilon=1,\\ \sinh\xi,&\varepsilon=-1,\end{cases}
\end{equation}
so that $\cc^2+\varepsilon\s^2=1$, $\s'=\cc$ and $\cc'=-\varepsilon\s$.

Let $F:\Sigma^{n-1}\to\M$ be an oriented immersed hypersurface with unit normal $\nu$ tangent to $\M$. We use the second fundamental form $h_{ij}=-\la\partial_i\partial_jF,\nu\ra_\varepsilon$, with principal curvatures $k_1,\dots,k_{n-1}$, and the normalised mean curvature
\[
H=\frac1{n-1}\sum_{i=1}^{n-1}k_i=\frac1{n-1}\,g^{ij}h_{ij}=\frac1{n-1}\,\mathrm{div}_\Sigma\nu .
\]
With this convention, a geodesic sphere of radius $\varrho$ has $H=\cc(\varrho)/\s(\varrho)>0$ with respect to the outward normal, and a horosphere has $H=1$, in every dimension. A family $F_t$ is a solution to the IMCF if $\partial_tF=H^{-1}\nu$, up to tangential diffeomorphisms. With this normalisation, as in \cite{LRS}, the soliton equation below reads $H\la X,\nu\ra=1$. With the unnormalised curvature $(n-1)H$, every speed $\lambda$ below would be divided by $n-1$.

\begin{definition}
Let $X$ be a Killing field on $\M$ generated by $B\in\mathfrak{so}(\varepsilon)$, the Lie algebra of $\la\cdot,\cdot\ra_\varepsilon$. A hypersurface $\Sigma$ with $H>0$ is an \emph{IMCF soliton with respect to $X$} if $\Sigma_t=\exp(tB)\Sigma$ is a solution to the IMCF.
\end{definition}

Since $\exp(tB)$ is an isometry that preserves $X$, both $H$ and $\la X,\nu\ra$ are preserved by the motion. Hence $\Sigma$ is a soliton if and only if
\begin{equation}\label{eq:soliton}
H\,\la X,\nu\ra=1 \quad\text{on }\Sigma.
\end{equation}

\begin{proposition}\label{prop:noclosed}
There is no closed IMCF soliton in $\M$.
\end{proposition}

\begin{proposition}\label{prop:cmc}
Let $\Sigma$ be a connected IMCF soliton in $\M$, with respect to an arbitrary Killing field $X$, and assume that its mean curvature is constant. Then $\varepsilon=-1$, and $\Sigma$ is an open subset of a horosphere. In particular, $H\equiv1$ and $\la X,\nu\ra\equiv1$ on $\Sigma$.
\end{proposition}

\begin{remark}
If $H\equiv H_0$ on an IMCF soliton, then $\la H_0^2X,\nu\ra=H$, so $\Sigma$ is also a soliton of the mean curvature flow, with respect to $H_0^2X$. Hence, for $n=3$, Proposition~\ref{prop:cmc} in $\Hy^3$ also follows from \cite[Theorem~3.4 and Remark~3.5]{LRS}. Our proof is direct, works in every dimension and also covers $\Sp^n$.
\end{remark}

\subsection{Rotation hypersurfaces}\label{subsec:classes}
Following do Carmo and Dajczer \cite{dCD83}, let $P^2\subset\E$ be a $2$-plane and let $G$ be the identity component of the group of isometries of $\M$ (restrictions of linear isometries of $\E$) that fix $P^2$ pointwise. A \emph{rotation hypersurface} is a $G$-invariant hypersurface. In $\Hy^n$, the rotation is of \emph{spherical}, \emph{hyperbolic} or \emph{parabolic type} according as $\la\cdot,\cdot\ra_\varepsilon$ restricted to $P^2$ is Lorentzian, Riemannian or degenerate; the names refer to the $G$-orbits. In $\Sp^n$ all planes are conjugate, and we also call this single class the spherical type. Up to an isometry, the three types are described as follows.

\smallskip
\emph{(S) Spherical type, $\varepsilon=\pm1$.} Let $e_\varepsilon(\theta)=(\cos\theta,\sin\theta)$ if $\varepsilon=1$ and $e_\varepsilon(\theta)=(\cosh\theta,\sinh\theta)$ if $\varepsilon=-1$, and let $\Sp^{n-2}\subset\R^{n-1}$ be the unit sphere. The map
\begin{equation}\label{eq:coordsS}
x(r,\theta,\omega)=\big(\cc(r)\,e_\varepsilon(\theta),\ \s(r)\,\omega\big),\qquad \omega\in\Sp^{n-2},
\end{equation}
gives coordinates on $\M$ with $r\in(0,\frac\pi2)$ if $\varepsilon=1$ and $r>0$ if $\varepsilon=-1$, in which
\[
g=dr^2+\cc(r)^2d\theta^2+\s(r)^2g_{\Sp^{n-2}} .
\]
Here $P^2=\mathrm{span}\{e_0,e_1\}$, $G=SO(n-1)$ acts on $\omega$, and its orbits are the round spheres $\{r,\theta\ \text{fixed}\}$. The \emph{axis} $\{r=0\}=\M\cap P^2$ is a geodesic. When $\varepsilon=1$, the \emph{far axis} $\{r=\frac\pi2\}=\Sp^n\cap(P^2)^\perp$ is a totally geodesic $\Sp^{n-2}$. We set $X_*=\partial_\theta$, the rotation of $P^2$ if $\varepsilon=1$ and the hyperbolic translation along the axis if $\varepsilon=-1$.

\smallskip
\emph{(H) Hyperbolic type, $\varepsilon=-1$.} Let $\Hy^{n-2}$ be the hyperbolic space in the first $n-1$ coordinates. The map
\begin{equation}\label{eq:coordsH}
x(r,p,\varphi)=\big(\cosh r\,p,\ \sinh r\,(\cos\varphi,\sin\varphi)\big),\qquad p\in\Hy^{n-2},
\end{equation}
gives coordinates on $\Hy^n$ with $r>0$, in which
\[
g=dr^2+\cosh^2r\,g_{\Hy^{n-2}}+\sinh^2r\,d\varphi^2 .
\]
Here $P^2=\mathrm{span}\{e_{n-1},e_n\}$, $G=SO^+(1,n-2)$ acts on $p$, and its orbits are the hyperbolic spaces $\{r,\varphi\ \text{fixed}\}$. The \emph{axis} $\{r=0\}=\Hy^n\cap(P^2)^\perp$ is a totally geodesic $\Hy^{n-2}$, of codimension two. We set $X_*=\partial_\varphi$, the rotation of $P^2$ around the axis.

\smallskip
\emph{(P) Parabolic type, $\varepsilon=-1$.} In the upper half-space model $\U^n=\{(z_1,\dots,z_{n-1},h):h>0\}$, with metric $h^{-2}(|dz|^2+dh^2)$, the group $G$ is the group of horizontal translations in $(z_2,\dots,z_{n-1})$, and its orbits are the horospherical planes $\{z_1,h\ \text{fixed}\}$. We set $X_*=\partial_{z_1}$.

\smallskip
For $\varepsilon=1$, the construction \eqref{eq:coordsH} with $p\in\Sp^{n-2}$ gives the spherical type again, after $r\mapsto\frac\pi2-r$. In each case $X_*$ commutes with $G$. If $Z$ is a Killing field tangent to the $G$-orbits, then $X=\lambda X_*+Z$ and $\lambda X_*$ give the same equation \eqref{eq:soliton} on $G$-invariant hypersurfaces. Replacing $\theta$ by $-\theta$ (or $\varphi$ by $-\varphi$, or $z_1$ by $-z_1$) if necessary, we assume from now on that
\[
X=\lambda X_*,\qquad \lambda>0 .
\]

In the types (S) and (H) we write $\vartheta$ for the angle $\theta$, respectively $\varphi$, and
\begin{equation}\label{eq:rhobeta}
(\rho,\beta)=(\cc,\s)\ \text{for (S)},\qquad (\rho,\beta)=(\s,\cc)\ \text{for (H)},
\end{equation}
so that $|X_*|=\rho(r)$, the orbit space carries the metric $dr^2+\rho(r)^2d\vartheta^2$, and the $G$-orbits have volume element proportional to $\beta(r)^{n-2}$. A $G$-invariant hypersurface is generated by a profile curve in the orbit space. We parametrise it by arc length $s$ and denote by $\sigma(s)$ the angle between the profile and $\partial_r$, with the orientation fixed by
\begin{equation}\label{eq:gaugen}
r'=\cos\sigma,\qquad \vartheta'=-\frac{\sin\sigma}{\rho(r)} .
\end{equation}

\begin{theorem}\label{thm:master}
Let $\Sigma$ be a rotation hypersurface of spherical or hyperbolic type. For the unit normal $N=\sin\sigma\,\partial_r+\rho^{-1}\cos\sigma\,\partial_\vartheta$,
\begin{equation}\label{eq:HXn}
(n-1)H=\sigma'+\big(\log(\rho\,\beta^{n-2})\big)'\sin\sigma,\qquad \la X,N\ra=\lambda\,\rho\cos\sigma .
\end{equation}
Hence $\Sigma$ is an IMCF soliton with respect to $X=\lambda X_*$ if, and only if,
\begin{equation}\label{eq:master}
\sigma'=D(r)\sin\sigma+\frac{\mathcal F(r)}{\cos\sigma},\qquad D=-\big(\log(\rho\,\beta^{n-2})\big)',\qquad \mathcal F=\frac{n-1}{\lambda\rho}.
\end{equation}
Equivalently, with $v=\sin\sigma$ and $r$ as parameter,
\begin{equation}\label{eq:scalar}
\frac{dv}{dr}=D(r)\,v+\frac{\mathcal F(r)}{\sqrt{1-v^2}} .
\end{equation}
Along every solution, the following flux identity holds:
\begin{equation}\label{eq:fluxid}
\big(\rho\,\beta^{n-2}v\big)'=\frac{(n-1)\,\beta^{n-2}}{\lambda\sqrt{1-v^2}} .
\end{equation}
\end{theorem}

Table~\ref{tab:classes} lists the three basic cases.

\begin{table}[h]
\centering
\renewcommand{\arraystretch}{1.3}
\begin{tabular}{lcc}
\toprule
class & drift $D(r)$ & forcing $\mathcal F(r)$\\
\midrule
spherical type in $\Sp^n$ (caps) & $\tan r-(n-2)\cot r$ & $\dfrac{n-1}{\lambda\cos r}$\\
spherical type in $\Hy^n$ (translators) & $-\big(\tanh r+(n-2)\coth r\big)$ & $\dfrac{n-1}{\lambda\cosh r}$\\
hyperbolic type in $\Hy^n$ (rotators) & $-\big(\coth r+(n-2)\tanh r\big)$ & $\dfrac{n-1}{\lambda\sinh r}$\\
\bottomrule
\end{tabular}
\caption{The equation \eqref{eq:scalar} in the basic classes.}\label{tab:classes}
\end{table}

The profile of a soliton can reach a point where $\cos\sigma=0$, that is, where $\la X,N\ra=0$. There $H\to\infty$. Our next result describes the hypersurface near such an orbit.

\begin{theorem}\label{thm:fronts}
Let $\sigma(s)\to\pi/2$ as $s\to s_*$ along a solution of \eqref{eq:master}, with $r(s)\to r_*$ in the interior, and let $\mathcal F_0=\mathcal F(r_*)$. Then $s_*<\infty$, $(n-1)H\sim\sqrt{\mathcal F_0/(2(s_*-s))}$, and the profile extends as a front whose Legendrian lift $(r,\vartheta,\sigma)$ is smooth and immersed. The profile has an ordinary $(2,3)$-cusp at the endpoint. Hence the hypersurface extends beyond the singular $G$-orbit as a front with a cuspidal edge along it. On each regular piece it satisfies \eqref{eq:soliton}. Moreover, the hypersurface is of class $C^{1,1/2}$ up to the singular orbit from each side, and not of class $C^2$.
\end{theorem}

We call such an orbit a \emph{singular orbit}: a \emph{singular sphere} $\Sp^{n-2}$ in the spherical type, a \emph{singular hyperbolic space} $\Hy^{n-2}$ in the hyperbolic type, and a \emph{singular horospherical plane} in the parabolic type. When $n=3$ they are circles, hypercycles and horocycles. As immersed hypersurfaces, solitons end at singular orbits. As fronts, they continue through them.

\subsection{Results in dimension $n$}
We now state the results for each type. Recall that $n\ge3$, $X=\lambda X_*$ and $\lambda>0$. We write
\begin{equation}\label{eq:kappan}
\kappa_n=\frac{2(n-1)}{(n-2)\,\lambda} .
\end{equation}

\begin{theorem}[Rotational translators in $\Hy^n$]\label{thm:transl}
Let $\varepsilon=-1$, and consider the spherical type. For each $\lambda>0$ there is a unique soliton $\Sigma_\lambda$ meeting the axis. It is smooth across the axis, and its mean curvature at the vertex is $1/\lambda$. Moreover:
\begin{enumerate}[label=(\roman*)]
\item if $\lambda>1$, then $\Sigma_\lambda$ is an entire bowl. Along the profile, $\sin\sigma=\kappa_ne^{-r}(1+o(1))$ and $H=\frac2\lambda e^{-r}(1+o(1))$. The hypersurface converges to the totally geodesic hyperplane $\{\theta=\theta_\infty\}$ orthogonal to the axis, with $\theta_\infty-\theta(r)=\kappa_ne^{-2r}(1+o(1))$;
\item if $\lambda=1$, then $\Sigma_\lambda$ is a horosphere;
\item if $\lambda<1$, then $\Sigma_\lambda$ reaches a singular sphere at a finite distance $r_+(\lambda)$ from the axis, $r_+$ is increasing in $\lambda$, and $r_+(\lambda)\to\infty$ as $\lambda\uparrow1$.
\end{enumerate}
By Proposition~\ref{prop:cmc}, the horosphere is the only one among these solitons with constant mean curvature.
\end{theorem}

The three cases are shown in Figure~\ref{fig:transl}.

In $\Hy^n$ the outer ends of both types share the following dichotomy. Let
\[
z_*(r)=\frac{\mathcal F(r)^2}{2(n-1)^2}=\frac{2}{\lambda^2}e^{-2r}\big(1+O(e^{-2r})\big).
\]

\begin{theorem}[Outer ends in $\Hy^n$]\label{thm:outer}
Let $\varepsilon=-1$, and consider the spherical or the hyperbolic type. Every solution of \eqref{eq:scalar} defined on $[R,\infty)$ is of exactly one of the following types:
\begin{enumerate}[label=(\Alph*)]
\item \emph{bowl end:} $\sin\sigma=\kappa_ne^{-r}+O(e^{-2r})$, $\vartheta\to\vartheta_\infty$ with error $O(e^{-2r})$, $H\to0$, and the end is complete and asymptotic to the totally geodesic hypersurface $\{\vartheta=\vartheta_\infty\}$, which is orthogonal to the axis in the spherical type and contains the axis in the hyperbolic type;
\item \emph{horospherical end:} $1-\sin\sigma=z_*(r)(1+o(1))$, the end is complete, $H\to1$ and $d\vartheta/dr\to-\lambda$.
\end{enumerate}
For $R$ large both types occur, and there is exactly one solution $v_B$ of type (B). It is a threshold: a solution with $v(R)<v_B(R)$ is of type (A), and a solution with $v(R)>v_B(R)$ reaches a singular orbit at a finite distance from the axis. For the rotational translators with $\lambda=1$, the solution of type (B) is the horosphere.
\end{theorem}

\begin{theorem}[Rotators of hyperbolic type in $\Hy^n$]\label{thm:hyp}
Let $\varepsilon=-1$, and consider the hyperbolic type. A solution of \eqref{eq:scalar} continued towards the axis is called an \emph{inward solution}.
\begin{enumerate}[label=(\roman*)]
\item \emph{Outer ends.} Theorem~\ref{thm:outer} applies, and both types of ends exist.
\item \emph{Inner ends, $\lambda>2(n-1)$.} Let $v_1<\frac1{\sqrt2}<v_2$ be the roots of $v\sqrt{1-v^2}=(n-1)/\lambda$, and let $\sigma_i=\arcsin v_i$, so that $\sin2\sigma_i=2(n-1)/\lambda$ and $\sigma_1+\sigma_2=\frac\pi2$. There is exactly one inward solution $v_e$, defined near the axis, with $v_e\to v_1$. Every inward solution with $v>v_e$ near the axis reaches the axis with $\sigma\to\sigma_2$, and every inward solution with $v<v_e$ reaches a singular orbit at a positive distance $r_-$ from the axis; beyond it, the front reaches the axis with $\sigma\to\sigma_2\pmod\pi$. In all cases the axis is reached after a finite arc length, the profile is a logarithmic spiral around the axis, of angle $\sigma_2$, or $\sigma_1$ for the exceptional solution $v_e$, and $H\sim1/(\lambda r\cos\sigma_i)\to\infty$.
\item \emph{Inner ends, $\lambda=2(n-1)$.} There is exactly one inward solution $v_c$, defined near the axis, with
\[
v_c-\frac1{\sqrt2}=-\frac{2n-3}{4\sqrt2}\,r^2+O(r^4).
\]
Every inward solution with $v>v_c$ near the axis reaches the axis after a finite arc length, with
\[
\sin\sigma-\frac1{\sqrt2}=\frac{1+o(1)}{2\sqrt2\,|\log r|},
\]
and every inward solution with $v<v_c$ reaches a singular orbit at some $r_->0$. In the first two cases $\sigma\to\frac\pi4$, and the profile is again a logarithmic spiral. In particular, every inward solution with $\frac1{\sqrt2}<\sin\sigma(r_1)<1$ at some $r_1>0$ reaches the axis.
\item \emph{Inner ends, $\lambda<2(n-1)$.} Every inward solution reaches a singular orbit at a positive distance $r_-$ from the axis. Continued as a front, it has infinitely many cuspidal edges, which accumulate at the axis; it winds infinitely many times around the axis, and its total length is finite.
\end{enumerate}
\end{theorem}

\begin{theorem}[Rotational solitons in $\Sp^n$]\label{thm:sphere}
Let $\varepsilon=1$. For each $\lambda>0$ there is a unique cap $\Sigma_\lambda$ meeting the axis. It is smooth across the axis, with vertex mean curvature $1/\lambda$. It reaches a singular sphere at some $r_+(\lambda)<\pi/2$. In particular, no rotational soliton of $\Sp^n$ containing a point of the axis is complete. Beyond the singular sphere, the front reaches the far axis. Moreover, with $x=\frac\pi2-r$:
\begin{enumerate}[label=(\roman*)]
\item if $\lambda<2(n-1)$, the angle $\sigma$ grows at least like $c\log(1/x)$, $c=\frac{n-1}\lambda-\frac12$. The front winds infinitely many times around the far axis and has infinitely many cuspidal edges, but its total length is finite;
\item if $\lambda>2(n-1)$, then $\sigma\to\sigma_*=\frac\pi2+\frac12\arcsin\frac{2(n-1)}\lambda \pmod \pi$, except possibly for one exceptional orbit of the far-axis dynamics, along which $\sigma\to\sigma_{**}=\pi-\frac12\arcsin\frac{2(n-1)}\lambda$. In both cases, the front arrives at the far axis along a logarithmic spiral, which winds infinitely many times around it at a constant angle;
\item if $\lambda=2(n-1)$, the front has finitely many cuspidal edges, and then $\sigma\to\frac{3\pi}4\pmod\pi$ with
\[
\frac{3\pi}4-\sigma=\frac{1+o(1)}{2|\log x|}\pmod\pi .
\]
\end{enumerate}
\end{theorem}

Figures~\ref{fig:boost} and~\ref{fig:sphere} illustrate Theorems~\ref{thm:hyp} and~\ref{thm:sphere} for $n=3$. The radius of the singular sphere of the caps has the following sharp asymptotics.

\begin{theorem}\label{thm:rplus}
Let $r_+(\lambda)$ be the distance from the axis to the singular sphere of the rotational cap $\Sigma_\lambda$ (Theorems~\ref{thm:transl}(iii) and~\ref{thm:sphere}).
\begin{enumerate}[label=(\roman*)]
\item As $\lambda\to0$, for both $\varepsilon=1$ and $\varepsilon=-1$,
\[
r_+(\lambda)=c_0(n)\lambda+O(\lambda^3).
\]
Here $\frac12<c_0(n)<1$ is the radius of the singular sphere of the rotational IMCF translator cap of $\R^n$ with speed $1$, and $r_+(\lambda)/\lambda$ is a smooth function of $\varepsilon\lambda^2$.
\item If $\varepsilon=1$, then, as $\lambda\to\infty$,
\[
\frac\pi2-r_+(\lambda)=\frac{K_n}\lambda-\frac{2(n-1)^2K_n\,\log\lambda}{\lambda^3}+O(\lambda^{-3}),\qquad K_n=(n-1)\int_0^{\pi/2}\cos^{n-2}t\,dt=\frac{\sqrt\pi\,\Gamma(\frac{n+1}2)}{\Gamma(\frac n2)} .
\]
\item In both space forms, the $(n-1)$-dimensional volume of the cap is
\[
\mathrm{Vol}(\Sigma_\lambda)=\frac{\omega_{n-2}}{n-1}\,\lambda\,\cc(r_+)\,\s(r_+)^{n-2},
\]
where $\omega_{n-2}$ is the volume of the unit sphere $\Sp^{n-2}$. If $\varepsilon=1$, then $\mathrm{Vol}(\Sigma_\lambda)\to\frac12\mathrm{Vol}(\Sp^{n-1})$ as $\lambda\to\infty$.
\end{enumerate}
\end{theorem}

For instance, $K_3=2$, $K_4=\frac{3\pi}4$, $K_5=\frac83$ and $K_6=\frac{15\pi}{16}$. By Stirling's formula, $K_n=\sqrt{\pi n/2}\,(1+O(1/n))$. As $\lambda\to\infty$ the caps of $\Sp^n$ converge to a totally geodesic hemisphere bounded by the far axis, and (ii) measures how close to the far axis the singular sphere is. There is no term of order $\lambda^{-2}$, and the first correction is logarithmic.

\begin{theorem}[Parabolic solitons in $\Hy^n$]\label{thm:parab}
Consider the parabolic type, and let the profile $(z_1(s),h(s))$ be parametrised by hyperbolic arc length, with $z_1'=h\cos\alpha$ and $h'=h\sin\alpha$. After a dilation we may take $\lambda=n-1$, and the soliton equation becomes
\begin{equation}\label{eq:parab}
\alpha'=-(n-1)\cos\alpha+\frac{h}{\sin\alpha},\qquad \alpha\in(0,\pi).
\end{equation}
Every maximal solution reaches a singular horospherical plane $\{h=h_{\max},\,z_1=z_1^*\}$ in forward time, at finite arc length. In backward time, exactly one of the following holds:
\begin{enumerate}[label=(\roman*)]
\item the profile reaches a singular horospherical plane, with $\alpha\to0$;
\item the end is complete, $h\to0$, $\alpha\to\pi/2$, $z_1\to z_1^\infty$, $H\to0$, and the hypersurface is asymptotic to the vertical totally geodesic hyperplane $\{z_1=z_1^\infty\}$;
\item the profile is the unique separatrix with $(h,\alpha)\to(0,0)$, along which $\alpha=\frac{h}{n-1}(1+o(1))$. This end is complete, $H\to1$, and it converges to the ideal point $\infty$ fixed by the parabolic group.
\end{enumerate}
Cases (i) and (ii) occur on open sets of initial data. Case (iii) is the common boundary between them.
\end{theorem}

The phase portrait of \eqref{eq:parab} and the three types of solitons are shown in Figure~\ref{fig:parab}.

\subsection{The three-dimensional case}\label{subsec:three}
Let $n=3$. In the coordinates
\begin{equation}\label{eq:coords}
x(r,\theta,\varphi)=\begin{cases}(\cos r\cos\theta,\ \cos r\sin\theta,\ \sin r\cos\varphi,\ \sin r\sin\varphi),& \varepsilon=1,\ r\in(0,\tfrac\pi2),\\
(\cosh r\cosh\theta,\ \cosh r\sinh\theta,\ \sinh r\cos\varphi,\ \sinh r\sin\varphi),&\varepsilon=-1,\ r>0,\end{cases}
\end{equation}
the metric is $g=dr^2+\cc(r)^2d\theta^2+\s(r)^2d\varphi^2$. These are the coordinates \eqref{eq:coordsS} with $\omega=(\cos\varphi,\sin\varphi)$, and, when $\varepsilon=-1$, also the coordinates \eqref{eq:coordsH} with $p=(\cosh\theta,\sinh\theta)$. The rotation group is now one-dimensional, and a rotation surface is invariant under one Killing field $Y$ commuting with $X$. The commuting pairs are classified as follows. This classification is classical: for $\varepsilon=1$ it is the conjugacy of maximal tori in $SO(4)$ (see e.g.\ \cite{BtD}), and for $\varepsilon=-1$ it follows from the classification of one-parameter groups of isometries of $\Hy^3$ (elliptic, hyperbolic, loxodromic and parabolic) and of their centralisers (see e.g.\ \cite{Be,Ra}). We include a short proof in Section~\ref{sec:three} only for completeness.

\begin{proposition}\label{prop:pairs}
Let $X,Y$ be linearly independent commuting Killing fields on $\mathbb M^3(\varepsilon)$.
\begin{enumerate}[label=(\roman*)]
\item If $\varepsilon=1$, or if $\varepsilon=-1$ and $\mathrm{span}\{X,Y\}$ contains a non-parabolic element, then, up to an isometry,
$Y=a\,\partial_\theta+b\,\partial_\varphi$ and $X=p\,\partial_\theta+q\,\partial_\varphi$.
\item Otherwise $\varepsilon=-1$ and, in the upper half-space model, up to an isometry, $Y=\partial_{z_2}$ and $X=\lambda\,\partial_{z_1}+\mu\,\partial_{z_2}$, $\lambda\neq0$.
\end{enumerate}
\end{proposition}

In case (i), the norm $|Y|$ depends only on $r$, and we write
\begin{equation}\label{eq:Delta}
|Y|^2=a^2\cc(r)^2+b^2\s(r)^2,\qquad \lambda=pb-qa,\qquad \psi=b\theta-a\varphi .
\end{equation}
The classes $(a,b)=(0,1)$ and $(a,b)=(1,0)$ are the spherical and the hyperbolic types for $n=3$, with $X_*=\partial_\theta$ and $X_*=-\partial_\varphi$, and case (ii) is the parabolic type. We parametrise the profile by arc length as before, with $r'=\cos\sigma$.

\begin{theorem}\label{thm:master3}
Let $n=3$ and let $\Sigma$ be a $Y$-invariant surface as in Proposition~\ref{prop:pairs}(i). For a suitable choice of unit normal $N$,
\begin{equation}\label{eq:HX}
2H=\sigma'+\Big(\frac{\cc}{\s}-\varepsilon\frac{\s}{\cc}\Big)\sin\sigma,\qquad
\la X,N\ra=\lambda\,\frac{\s\,\cc\,\cos\sigma}{|Y|}.
\end{equation}
Hence $\Sigma$ is an IMCF soliton with respect to $X$ if, and only if, $\lambda\neq0$ and
\begin{equation}\label{eq:master3}
\sigma'=\Big(\varepsilon\frac{\s}{\cc}-\frac{\cc}{\s}\Big)\sin\sigma+\frac{2|Y|}{\lambda\,\s\,\cc\,\cos\sigma}.
\end{equation}
\end{theorem}

The mean curvature in \eqref{eq:HX} does not depend on $(a,b)$. Only the forcing term in \eqref{eq:master3} sees the pitch, through $|Y|^2$, and only the combination $\lambda=pb-qa$ of the soliton field enters. For $(a,b)=(0,1)$ and $(1,0)$, \eqref{eq:master3} is \eqref{eq:master} with $n=3$. The remaining basic class is the Hopf class $a=b=1$ in $\Sp^3$, with drift $-2\cot(2r)$ and forcing $4/(\lambda\sin 2r)$.

The general theorems give the following three-dimensional statements, with $\kappa_3=4/\lambda$, threshold $2(n-1)=4$, $K_3=2$ and parabolic normalisation $\lambda=2$.

\begin{corollary}[Rotational translators in $\Hy^3$]\label{cor:transl3}
Let $\varepsilon=-1$ and $(a,b)=(0,1)$. For each $\lambda>0$ there is a unique soliton $\Sigma_\lambda$ meeting the axis, smooth, with vertex mean curvature $1/\lambda$. If $\lambda>1$, it is an entire bowl with $\sin\sigma=\frac4\lambda e^{-r}(1+o(1))$, $H=\frac2\lambda e^{-r}(1+o(1))$ and $\theta_\infty-\theta(r)=\frac4\lambda e^{-2r}(1+o(1))$. If $\lambda=1$, it is a horosphere. If $\lambda<1$, it reaches a singular circle at $r_+(\lambda)$, which increases to $\infty$ as $\lambda\uparrow1$.
\end{corollary}

\begin{corollary}[Outer ends in $\Hy^3$]\label{cor:outer3}
Let $\varepsilon=-1$ and $(a,b)\neq(0,0)$, and let $\kappa=4\sqrt{a^2+b^2}/\lambda$. Every solution of \eqref{eq:master3} defined on $[R,\infty)$ either is a bowl end, with $\sin\sigma=\kappa e^{-r}+O(e^{-2r})$, $\psi\to\psi_\infty$ with error $O(e^{-2r})$ and $H\to0$, asymptotic to the minimal surface $\{\psi=\psi_\infty\}$, or is a horospherical end, with $1-\sin\sigma=\frac{\kappa^2}8e^{-2r}(1+o(1))$, $H\to1$ and $d\psi/dr\to-\lambda$. For $R$ large there is exactly one horospherical end $v_B$; the solutions with $v(R)<v_B(R)$ are bowl ends, and those with $v(R)>v_B(R)$ reach a singular orbit.
\end{corollary}

\begin{corollary}[Boost-invariant rotators in $\Hy^3$]\label{cor:boost3}
Let $\varepsilon=-1$ and $(a,b)=(1,0)$. The outer ends are as in Corollary~\ref{cor:outer3}, the bowl ends converging to totally geodesic planes that contain the axis. Let $v_1<v_2$ be the roots of $v\sqrt{1-v^2}=2/\lambda$ when $\lambda>4$. If $\lambda>4$, there is exactly one inward solution $v_e$ with $v_e\to v_1$; the inward solutions with $v>v_e$ near the axis reach the axis along logarithmic spirals with $\sin\sigma\to v_2$, and those with $v<v_e$ reach a singular hypercycle, beyond which the front reaches the axis in the same way. If $\lambda=4$, there is exactly one inward solution $v_c$ with $v_c-\frac1{\sqrt2}=-\frac{3}{4\sqrt2}r^2+O(r^4)$; the inward solutions with $v>v_c$ reach the axis with $\sin\sigma-\frac1{\sqrt2}=\frac{1+o(1)}{2\sqrt2|\log r|}$, and those with $v<v_c$ reach a singular hypercycle. If $\lambda<4$, every inward solution reaches a singular hypercycle at a positive distance from the axis, and the front continues through infinitely many cuspidal edges accumulating at the axis.
\end{corollary}

\begin{corollary}[Rotational solitons in $\Sp^3$]\label{cor:sphere3}
Let $\varepsilon=1$ and $(a,b)=(0,1)$. The cap $\Sigma_\lambda$ reaches a singular circle at $r_+(\lambda)<\frac\pi2$, and the front reaches the far axis. If $\lambda<4$, it winds infinitely many times around the far axis with infinitely many cuspidal edges and finite length, and $\sigma$ grows at least like $(\frac2\lambda-\frac12)\log(1/x)$. If $\lambda>4$, it arrives along a logarithmic spiral with $\sigma\to\frac\pi2+\frac12\arcsin\frac4\lambda\pmod\pi$, except possibly for one orbit with $\sigma\to\pi-\frac12\arcsin\frac4\lambda$. If $\lambda=4$, after finitely many cuspidal edges, $\frac{3\pi}4-\sigma=\frac{1+o(1)}{2|\log x|}\pmod\pi$.
\end{corollary}

\begin{corollary}\label{cor:rplus3}
For $n=3$, $r_+(\lambda)=c_0\lambda+O(\lambda^3)$ as $\lambda\to0$, where $c_0=c_0(3)\in(\frac12,1)$ is the radius of the singular circle of the rotational IMCF translator cap of $\R^3$ with speed $1$, and, in $\Sp^3$, $\frac\pi2-r_+(\lambda)=\frac2\lambda-\frac{16\log\lambda}{\lambda^3}+O(\lambda^{-3})$ as $\lambda\to\infty$. The area of the cap is $\pi\lambda\,\cc(r_+)\,\s(r_+)$, and it tends to $2\pi$ in $\Sp^3$ as $\lambda\to\infty$.
\end{corollary}

\begin{corollary}[Parabolic solitons in $\Hy^3$]\label{cor:parab3}
In case (ii) of Proposition~\ref{prop:pairs}, after a dilation $\lambda=2$, the soliton equation is $\alpha'=-2\cos\alpha+h/\sin\alpha$, and the conclusions of Theorem~\ref{thm:parab} hold, with singular horocycles and separatrix $\alpha=\frac h2(1+o(1))$.
\end{corollary}

The helicoidal and Hopf classes have no analogue in Theorems~\ref{thm:transl}--\ref{thm:parab}.

\begin{theorem}[Helicoidal solitons]\label{thm:helic}
Let $n=3$ and $ab\neq0$.
\begin{enumerate}[label=(\roman*)]
\item If $\varepsilon=-1$, every end with $r\to\infty$ either satisfies $\sin\sigma=\frac{4\sqrt{a^2+b^2}}{\lambda}e^{-r}(1+o(1))$ and is asymptotic, at the rate $O(e^{-2r})$, to the minimal helicoid of pitch $(a,b)$, or is the unique horospherical end of Corollary~\ref{cor:outer3}. Near the axis, the conclusions of Corollary~\ref{cor:boost3} hold with $\lambda/|a|$ in place of $\lambda$, with $v_c-\frac1{\sqrt2}=O(r^2)$ in the borderline case $\lambda=4|a|$, and with singular orbits of $Y$ (helices) in place of hypercycles.
\item If $\varepsilon=1$, no helicoidal soliton is complete. Each end reaches a singular orbit of $Y$ at finite arc length, or, when $\lambda\ge4|a|$ (respectively $\lambda\ge4|b|$), it may reach the axis (respectively the far axis) along a logarithmic spiral.
\end{enumerate}
\end{theorem}

\begin{theorem}[Hopf solitons and the inverse curve shortening flow]\label{thm:hopf}
Let $\varepsilon=1$, $n=3$ and $a=b=1$, so that $Y$ is the Hopf field. A Hopf-invariant surface is an IMCF soliton with respect to $X=p\partial_\theta+q\partial_\varphi$ if, and only if, its projection $\gamma\subset\Sp^2(\frac12)$ is a rotator of the inverse curve shortening flow $\partial_t\gamma=\kappa^{-1}n$ with angular speed $\omega=\frac12(p-q)=\frac\lambda2$. There is no complete rotator. Every rotator reaches, in each direction, either a point where $\kappa\to\infty$ (the image of a singular circle), or, only if $\omega\ge2$, one of the two fixed points of the rotation, along a logarithmic spiral.
\end{theorem}

Examples of helicoidal and Hopf solitons are shown in Figures~\ref{fig:helic} and~\ref{fig:hopf}.

\begin{remark}\label{rem:rank}
The helicoidal and Hopf classes use that, for $n=3$, the rotation group and the soliton field generate an abelian group of isometries. For $n\ge4$, an abelian group of isometries of $\Sp^n$ with orbits of codimension two, together with a commuting Killing field, would have rank at least $n-1$, while the rank of $SO(n+1)$ is $\lfloor\frac{n+1}2\rfloor<n-1$. So no such class exists in $\Sp^n$ for $n\ge4$. In $\Hy^n$ the parabolic type is abelian in every dimension.
\end{remark}

\section{Preliminaries and general properties}\label{sec:prelim}

Since $\E$ is flat and $\la x,\nu\ra=0$, the Levi-Civita connection of $\M$ is $\overline\nabla_UV=\mathrm D_UV+\varepsilon\la U,V\ra x$, where $\mathrm D$ is the flat connection of $\E$, and the second fundamental form of $\Sigma\subset\M$ is
\[
h_{ij}=-\la\partial_i\partial_jF,\nu\ra=\la F_j,\overline\nabla_{\partial_i}\nu\ra .
\]
In particular $\mathrm D_U\nu=\overline\nabla_U\nu=A(U)$ for $U$ tangent to $\Sigma$, where $A$ is the shape operator.

\begin{lemma}\label{lem:div}
Let $X$ be a Killing field on $\M$, and let $\Sigma\subset\M$ be a hypersurface. Write $X=X^\top+\la X,\nu\ra\nu$ along $\Sigma$. Then $\mathrm{div}_\Sigma X^\top=-(n-1)H\la X,\nu\ra$.
\end{lemma}

\begin{proof}
Let $\{e_1,\dots,e_{n-1}\}$ be a local orthonormal frame on $\Sigma$. Since $\overline\nabla X$ is skew-adjoint, $\sum_i\la\overline\nabla_{e_i}X,e_i\ra=0$. On the other hand, $\sum_i\la\overline\nabla_{e_i}(\la X,\nu\ra\nu),e_i\ra=\la X,\nu\ra\,\mathrm{div}_\Sigma\nu=(n-1)H\la X,\nu\ra$. Subtracting the two identities proves the lemma.
\end{proof}

\begin{proof}[Proof of Proposition~\ref{prop:noclosed}]
On a soliton, Lemma~\ref{lem:div} and \eqref{eq:soliton} give $\mathrm{div}_\Sigma X^\top=-(n-1)$. If $\Sigma$ were closed, the divergence theorem would give $0=-(n-1)\,\mathrm{Vol}(\Sigma)$, which is impossible.
\end{proof}

\begin{remark}
For mean curvature flow solitons, $(n-1)H=\la X,\nu\ra$ in the unnormalised convention, and the same argument gives $\int_\Sigma H^2=0$, so closed solitons are minimal \cite{MMRS}. The IMCF normalisation excludes closed examples altogether.
\end{remark}

The next lemma describes parallel hypersurfaces; compare \cite{RT}, where parallel families are used for the mean curvature flow, and \cite{AT} for the IMCF.

\begin{lemma}\label{lem:parallel}
Let $p\in\Sigma$, and let $e_1,\dots,e_{n-1}$ be principal directions at $p$, with principal curvatures $k_1,\dots,k_{n-1}$. For small $d$, the map $F_d=\cc(d)F+\s(d)\nu$ is an immersion near $p$, with unit normal $\nu_d=-\varepsilon\s(d)F+\cc(d)\nu$, and its principal curvatures at $p$, with respect to $\nu_d$, are
\[
k_i(d)=\frac{\cc(d)\,k_i-\varepsilon\s(d)}{\cc(d)+\s(d)\,k_i},\qquad\text{with}\qquad \frac{d}{dd}k_i(d)=-\big(k_i(d)^2+\varepsilon\big).
\]
\end{lemma}

\begin{proof}
The curve $d\mapsto F_d(p)$ is the normal geodesic through $F(p)$. Since $\mathrm D_{e_i}\nu=k_ie_i$, we get $dF_d(e_i)=(\cc+\s k_i)e_i$ and $d\nu_d(e_i)=(-\varepsilon\s+\cc k_i)e_i$. One checks directly that $\la\nu_d,\nu_d\ra=\cc^2+\varepsilon\s^2=1$, $\la\nu_d,F_d\ra=0$ and $\la\nu_d,e_i\ra=0$. Hence $\nu_d$ is a unit normal, and the shape operator of $F_d$ is diagonal in the basis $e_i$, with eigenvalues $(\cc k_i-\varepsilon\s)/(\cc+\s k_i)$. The derivative formula is a direct computation using $\s'=\cc$ and $\cc'=-\varepsilon\s$.
\end{proof}

\begin{proof}[Proof of Proposition~\ref{prop:cmc}]
Let $H\equiv H_0$, so that $f=\la X,\nu\ra=1/H_0$ is constant. The family $\Sigma_t=\exp(tB)\Sigma$ solves the IMCF, and its normal speed $\la X,\nu_t\ra$ equals the constant $f$, because $\exp(tB)$ is an isometry preserving $X$. After a tangential reparametrisation $\tilde F_t$, we have $\partial_t\tilde F=f\nu_t$. For a normal variation, $\overline\nabla_t\nu_t=-\nabla f=0$ \cite{Hu86,HP}, hence $\overline\nabla_t\partial_t\tilde F=0$ and the curves $t\mapsto\tilde F_t(p)$ are the normal geodesics. Thus $\tilde F_t=F_{ft}$ is the parallel family of Lemma~\ref{lem:parallel}. Since $\tilde F_t(\Sigma)=\exp(tB)\Sigma$, its mean curvature is $H_0$ at every point. By Lemma~\ref{lem:parallel}, at each $p\in\Sigma$,
\[
0=\frac{d}{dd}\sum_ik_i(d)=-\sum_i\big(k_i(d)^2+\varepsilon\big)\qquad\text{for all small }d .
\]
If $\varepsilon=1$ this is impossible. Let $\varepsilon=-1$, and write $\tau=\tanh d$. Then $k_i(d)=(k_i+\tau)/(1+k_i\tau)$ and $1-k_i(d)^2=(1-\tau^2)(1-k_i^2)/(1+k_i\tau)^2$, so
\[
\sum_{i=1}^{n-1}\frac{1-k_i^2}{(1+k_i\tau)^2}=0\qquad\text{for all small }\tau .
\]
This is a rational function of $\tau$, so it vanishes identically. Group the indices by the distinct values $\kappa$ of the $k_i$. A value $\kappa\notin\{0,1,-1\}$ produces a double pole at $\tau=-1/\kappa$ with nonzero coefficient, and distinct values give distinct poles; hence no such value occurs. A value $\kappa=0$ contributes a positive constant, while the values $\pm1$ contribute zero; hence $\kappa=0$ does not occur either. So every principal curvature of $\Sigma$ is $\pm1$. The multiplicities are then locally constant, and $\Sigma$ is isoparametric. If both values $1$ and $-1$ occurred, Cartan's identity for isoparametric hypersurfaces with two distinct principal curvatures in $\Hy^n$, $\kappa_1\kappa_2=1$ \cite{Ca}, would give $-1=1$. Hence $\Sigma$ is totally umbilic with $A=\pm\mathrm{Id}$. Since $H_0>0$, $A=\mathrm{Id}$, and $\Sigma$ is an open subset of a horosphere. Finally $H_0=1$ and $\la X,\nu\ra=1/H_0=1$.
\end{proof}

\section{The reduction in dimension $n$}\label{sec:master}

Consider a metric of the form
\begin{equation}\label{eq:warped}
g=dr^2+\rho(r)^2d\vartheta^2+\beta(r)^2g_{\mathcal O},
\end{equation}
where $(\mathcal O,g_{\mathcal O})$ is a Riemannian manifold of dimension $n-2$. In the coordinates \eqref{eq:coordsS} and \eqref{eq:coordsH}, the metric of $\M$ has this form, with $\mathcal O=\Sp^{n-2}$ or $\Hy^{n-2}$ and $(\rho,\beta)$ as in \eqref{eq:rhobeta}. The Koszul formula gives, for vector fields $V,W$ lifted from $\mathcal O$,
\begin{equation}\label{eq:christn}
\overline\nabla_{\partial_r}\partial_r=0,\quad \overline\nabla_{\partial_\vartheta}\partial_\vartheta=-\rho\rho'\partial_r,\quad \overline\nabla_{\partial_\vartheta}\partial_r=\frac{\rho'}{\rho}\partial_\vartheta,\quad \overline\nabla_V\partial_r=\frac{\beta'}{\beta}V,\quad \overline\nabla_V\partial_\vartheta=0,
\end{equation}
and $\overline\nabla_VW=\nabla^{\mathcal O}_VW-\beta\beta'g_{\mathcal O}(V,W)\partial_r$.

\begin{proof}[Proof of Theorem~\ref{thm:master}]
Let $e_r=\partial_r$ and $e_\vartheta=\rho^{-1}\partial_\vartheta$. By \eqref{eq:christn}, $\overline\nabla_{e_r}e_r=\overline\nabla_{e_r}e_\vartheta=0$, $\overline\nabla_{e_\vartheta}e_r=\frac{\rho'}{\rho}e_\vartheta$ and $\overline\nabla_{e_\vartheta}e_\vartheta=-\frac{\rho'}{\rho}e_r$. By \eqref{eq:gaugen}, the unit tangent of the profile is $T=\cos\sigma\,e_r-\sin\sigma\,e_\vartheta$, and $N=\sin\sigma\,e_r+\cos\sigma\,e_\vartheta$ is a unit normal of $\Sigma$, since it is orthogonal to $T$ and to the $G$-orbits.

\emph{Profile direction.} Since $\overline\nabla_Te_r=-\sin\sigma\frac{\rho'}{\rho}e_\vartheta$ and $\overline\nabla_Te_\vartheta=\sin\sigma\frac{\rho'}{\rho}e_r$,
\[
\overline\nabla_TT=-\sigma'N-\frac{\rho'}{\rho}\sin\sigma\big(\sin\sigma\,e_r+\cos\sigma\,e_\vartheta\big)=-\Big(\sigma'+\frac{\rho'}{\rho}\sin\sigma\Big)N .
\]
Hence $h(T,T)=-\la\overline\nabla_TT,N\ra=\sigma'+\frac{\rho'}{\rho}\sin\sigma$.

\emph{Orbit directions.} For $V$ tangent to the orbit, $h(V,V)=-\la\overline\nabla_VV,N\ra=\beta\beta'g_{\mathcal O}(V,V)\sin\sigma=\frac{\beta'}{\beta}\sin\sigma\,|V|^2$. Since also $\overline\nabla_TV$ and $\overline\nabla_VT$ have no component that could make $h(T,V)\neq0$, the orbit directions are principal, with principal curvature $\frac{\beta'}{\beta}\sin\sigma$ and multiplicity $n-2$.

\emph{Sum.} Adding, $(n-1)H=\sigma'+\big(\frac{\rho'}{\rho}+(n-2)\frac{\beta'}{\beta}\big)\sin\sigma$, which is \eqref{eq:HXn}. Since $|X_*|=\rho$, $\la X,N\ra=\lambda\rho\la e_\vartheta,N\ra=\lambda\rho\cos\sigma$. Substituting into \eqref{eq:soliton} gives \eqref{eq:master}, and $dv/dr=\cos\sigma\,\sigma'/\cos\sigma=\sigma'$ gives \eqref{eq:scalar}. Finally, $(\rho\beta^{n-2})'=-D\rho\beta^{n-2}$, so
\[
\big(\rho\beta^{n-2}v\big)'=-D\rho\beta^{n-2}v+\rho\beta^{n-2}\Big(Dv+\frac{n-1}{\lambda\rho\sqrt{1-v^2}}\Big)=\frac{(n-1)\beta^{n-2}}{\lambda\sqrt{1-v^2}} .\qedhere
\]
\end{proof}

\begin{remark}
Formulas \eqref{eq:HXn} were also checked independently, for $n=3,\dots,6$ and both signs of $\varepsilon$, by computing the unit normal and the second fundamental form directly from the embedding in $\E$, for random profiles. The two computations agree to the accuracy of the finite differences. The flux identity \eqref{eq:fluxid} is the integrated form of Lemma~\ref{lem:div}: $\lambda\rho\beta^{n-2}\sin\sigma$ is, up to the volume of $\mathcal O$, the flux of $X$ through an orbit, and $\int\beta^{n-2}ds$ is the volume swept by the profile.
\end{remark}

\begin{proof}[Proof of Theorem~\ref{thm:fronts}]
Let $w=\frac\pi2-\sigma$ and multiply \eqref{eq:master} by $\cos\sigma$. Since the drift term contributes $O(\sin\sigma\cos\sigma)=O(w)$, we get $ww'=-\mathcal F_0+O(w)+O(s_*-s)$. Hence $w=\sqrt{2\mathcal F_0(s_*-s)}(1+o(1))$, and, by \eqref{eq:HXn}, $(n-1)H=\sigma'+O(1)=\sqrt{\mathcal F_0/(2(s_*-s))}(1+o(1))$. Since $\sigma'$ is integrable, $s_*<\infty$ and the tangent has a limit. Since $\sigma-\frac\pi2=O((s_*-s)^{1/2})$, the profile is $C^{1,1/2}$ and not $C^2$.

To continue the profile, introduce the parameter $t$ with $ds=\cos\sigma\,dt$. Along the profile,
\begin{equation}\label{eq:systemt}
\frac{dr}{dt}=\cos^2\sigma,\qquad \frac{d\vartheta}{dt}=-\frac{\sin\sigma\cos\sigma}{\rho},\qquad
\frac{d\sigma}{dt}=\mathcal F(r)+D(r)\sin\sigma\cos\sigma .
\end{equation}
This vector field is smooth near $\sigma=\pi/2$, and $d\sigma/dt=\mathcal F_0>0$ there. Hence $(r,\vartheta,\sigma)$ is a smooth immersed Legendrian curve in the unit tangent bundle of the orbit space, so the profile is a front. Writing $\sigma=\frac\pi2+\mathcal F_0t+O(t^2)$, we obtain
\[
r-r_*=\tfrac13\mathcal F_0^2\,t^3+O(t^4),\qquad \vartheta-\vartheta_*=\frac{\mathcal F_0}{2\rho(r_*)}\,t^2+O(t^3),
\]
which is an ordinary $(2,3)$-cusp. The hypersurface is locally the product of this cusp with an open subset of the singular orbit, so it has a cuspidal edge. On the branch $\cos\sigma<0$, both $H$ and $\la X,N\ra$ change sign with respect to the continuous front normal, so \eqref{eq:soliton} still holds.
\end{proof}

\section{Spherical and hyperbolic types}\label{sec:rot}

\subsection{Caps}
Consider the spherical type. Then $\rho=\cc$, $\beta=\s$, and near $r=0$ we have $D=\varepsilon\frac\s\cc-(n-2)\frac\cc\s=-\frac{n-2}r+\frac{(n+1)\varepsilon}{3}r+O(r^3)$. Hence \eqref{eq:scalar} becomes
\begin{equation}\label{eq:frob}
\frac{dv}{dr}=\Big(-\frac{n-2}r+\frac{(n+1)\varepsilon}3r+O(r^3)\Big)v+\frac{n-1}\lambda\big(1+\tfrac\varepsilon2r^2+O(r^4)\big)\big(1+\tfrac12v^2+O(v^4)\big).
\end{equation}

\begin{lemma}\label{lem:cap}
Equation \eqref{eq:frob} has a unique solution with $v=O(r)$ as $r\to0$. It is odd, with $v=\frac{r}{\lambda}+O(r^3)$, and the corresponding hypersurface is smooth across the axis, with mean curvature $1/\lambda$ at the vertex. Every other solution behaves like $Cr^{-(n-2)}$, $C\neq0$, and reaches $|v|=1$ before the axis.
\end{lemma}

\begin{proof}
The homogeneous part $v'=-(n-2)v/r$ has exponent $-(n-2)<0$. Writing \eqref{eq:frob} as $(r^{n-2}v)'=r^{n-2}\Phi(r,v)$, with $\Phi(0,0)=\frac{n-1}\lambda$, the solutions with $v=O(r)$ are exactly the fixed points of
\[
\mathcal T v(r)=\frac1{r^{n-2}}\int_0^rt^{n-2}\,\Phi(t,v(t))\,dt
\]
in $\{v\in C^0[0,\delta]:|v(r)|\le 2r/\lambda\}$, and $\mathcal T$ is a contraction there for $\delta$ small. Since $\int_0^rt^{n-2}dt=r^{n-1}/(n-1)$, the fixed point satisfies $v=\frac r\lambda+O(r^3)$. Equation \eqref{eq:frob} is invariant under $(r,v)\mapsto(-r,-v)$, so this solution is odd.

The profile is $\theta=\theta(r)$, with $d\theta/dr=-v/(\cc\sqrt{1-v^2})$ by \eqref{eq:gaugen}. Hence $\theta-\theta_0=-\frac{r^2}{2\lambda}+O(r^4)$ is a smooth function of $r^2$, which is a smooth function on $\M$ near the axis. So the hypersurface is the regular level set $\{\theta-\theta_0-\mathfrak h(r^2)=0\}$, and it is smooth across the axis. At the vertex, \eqref{eq:HXn} gives $(n-1)H=v'-Dv\to\frac1\lambda+\frac{n-2}\lambda$, and $\la X,N\ra=\lambda\cc\cos\sigma\to\lambda$.
\end{proof}

The next identity is the integrated form of \eqref{eq:fluxid} along the cap.

\begin{lemma}[Flux identity]\label{lem:flux}
Along the cap of Lemma~\ref{lem:cap}, for $r<r_+$,
\begin{equation}\label{eq:flux}
\lambda\,\cc(r)\,\s(r)^{n-2}\sin\sigma(r)=(n-1)\int_0^r\frac{\s(\varrho)^{n-2}}{\sqrt{1-v(\varrho)^2}}\,d\varrho=(n-1)\int_0^{s(r)}\s(r(u))^{n-2}\,du .
\end{equation}
\end{lemma}

\begin{proof}
Integrate \eqref{eq:fluxid}, using that $\cc\s^{n-2}v\to0$ at $r=0$ and $ds=dr/\sqrt{1-v^2}$.
\end{proof}

Since the orbit through a point of the profile has volume $\omega_{n-2}\s(r)^{n-2}$, the right-hand side of \eqref{eq:flux} is $(n-1)/\omega_{n-2}$ times the volume of the part of the cap inside the tube of radius $r$ around the axis. Letting $r\to r_+$, where $\sin\sigma=1$, we obtain Theorem~\ref{thm:rplus}(iii):
\begin{equation}\label{eq:volcap}
\mathrm{Vol}(\Sigma_\lambda)=\frac{\omega_{n-2}}{n-1}\,\lambda\,\cc(r_+)\,\s(r_+)^{n-2}.
\end{equation}
This is Lemma~\ref{lem:div} integrated over the cap: its volume is $\frac1{n-1}$ times the flux of $X$ through the singular sphere.

We will use two elementary lemmas on asymptotically autonomous equations in the variable $\varsigma=-\log r$ (or $-\log x$).

\begin{lemma}\label{lem:atr}
Let $k,C>0$ and let $w$ solve $dw/d\varsigma=-kw+R(\varsigma,w)$ with $|R|\le C(w^2+e^{-2\varsigma})$ for $|w|\le\eta_0$. There are $\eta_1>0$ and $\varsigma_*$ such that, if $|w(\varsigma_1)|\le\eta_1/2$ and $\varsigma_1\ge\varsigma_*$, then $w$ is defined on $[\varsigma_1,\infty)$ and $w\to0$.
\end{lemma}

\begin{proof}
Let $\eta_1=\min\{\eta_0,k/4C\}$. While $|w|\le\eta_1$, Young's inequality gives $(w^2)'\le-kw^2+\frac{2C^2}{k}e^{-4\varsigma}$. By Gr\"onwall, $w^2\le e^{-k(\varsigma-\varsigma_1)}w(\varsigma_1)^2+\frac{C^2}{2k}e^{-4\varsigma_1}$. For $\varsigma_1$ large this stays below $\eta_1^2$, and both terms tend to $0$.
\end{proof}

The borderline cases are governed by the following elementary lemma, a non-autonomous saddle-node.

\begin{lemma}\label{lem:sn}
Let $a,C,\eta>0$ with $C\eta\le a/2$, and let $w$ solve
\[
\frac{dw}{d\varsigma}=-a\,w^2+R(\varsigma,w),\qquad |R(\varsigma,w)|\le C\big(|w|^3+e^{-2\varsigma}\big)\ \text{ for }|w|\le\eta .
\]
For every $\alpha\in(0,1/a)$ there is $\varsigma_\alpha$ with the following property. If $\varsigma_1\ge\varsigma_\alpha$ and $\alpha/\varsigma_1\le w(\varsigma_1)\le\eta$, then $w$ is defined on $[\varsigma_1,\infty)$, $\alpha/\varsigma\le w(\varsigma)\le\eta$, and $\varsigma\,w(\varsigma)\to1/a$.
\end{lemma}

\begin{proof}
Choose $\varsigma_\alpha$ so that, for $\varsigma\ge\varsigma_\alpha$,
\[
\alpha(1-a\alpha)\varsigma^{-2}>C\big(\alpha^3\varsigma^{-3}+e^{-2\varsigma}\big)\quad\text{and}\quad a\eta^2>C(\eta^3+e^{-2\varsigma}).
\]
Then $\alpha/\varsigma$ is a strict subsolution and $\eta$ is a strict supersolution, so $\alpha/\varsigma\le w\le\eta$. Next, for every $\eta'\in(0,\eta]$, the inequality $w'\le-(a-C\eta')w^2+Ce^{-2\varsigma}$ holds while $w\le\eta'$. Since $w\ge\alpha/\varsigma$, we have $(1/w)'\ge a-C\eta'-C\alpha^{-2}\varsigma^2e^{-2\varsigma}$. The last term is integrable, so $1/w\to\infty$, and hence $w\le\eta'$ eventually and $\limsup\varsigma w\le1/(a-C\eta')$. Similarly, $(1/w)'\le a+Cw+C\alpha^{-2}\varsigma^2e^{-2\varsigma}$ gives $\liminf\varsigma w\ge1/a$. Letting $\eta'\to0$ we conclude.
\end{proof}

\begin{figure}[tp]
\centering
\includegraphics[width=0.74\textwidth]{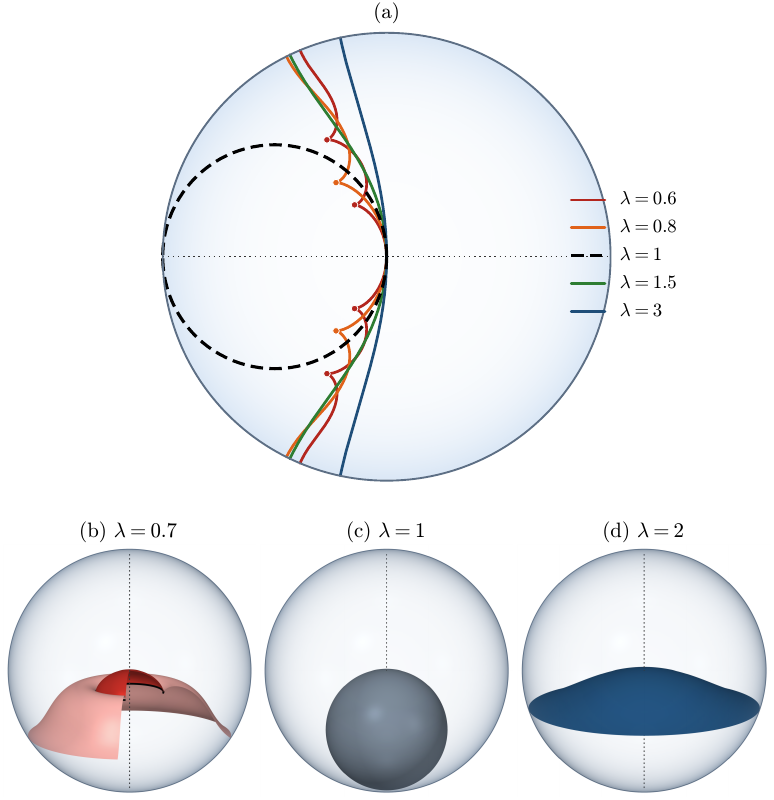}
\caption{Rotational translators of $\Hy^3$ (Theorem~\ref{thm:transl} and Corollary~\ref{cor:transl3}) in the Poincar\'e ball. (a) Meridian profiles of the caps $\Sigma_\lambda$; the axis is the dotted horizontal diameter. For $\lambda<1$ the profile reaches a singular circle (dots) and continues as a front; for $\lambda=1$ it is the horosphere (dashed), tangent to the ideal boundary; for $\lambda>1$ it is an entire bowl, which meets the ideal boundary along a circle. (b) The cap with $\lambda=0.7$ (dark) and its continuation as a front (light), cut open, with the singular circle. (c) The horosphere, $\lambda=1$. (d) The bowl with $\lambda=2$.}
\label{fig:transl}
\end{figure}

\subsection{Proof of Theorem~\ref{thm:transl}}
Let $\varepsilon=-1$ and write $F_\lambda(r,v)=-\big(\tanh r+(n-2)\coth r\big)v+\frac{n-1}{\lambda\cosh r\sqrt{1-v^2}}$. Then $\partial_\lambda F_\lambda<0$. Let $v_\lambda$ be the cap given by Lemma~\ref{lem:cap}. We will use that, for $r>0$,
\begin{equation}\label{eq:driftbound}
\tanh r+(n-2)\coth r>n-1,
\end{equation}
which holds because $(n-2)(\coth r-1)\ge\frac{2}{e^{2r}-1}>\frac{2}{e^{2r}+1}=1-\tanh r$.

\emph{The horosphere.} For $\lambda=1$, the function $v=\tanh r$ solves \eqref{eq:scalar} in every dimension. Indeed, $\sqrt{1-v^2}=\mathrm{sech}\,r$ and
\[
-\big(\tanh r+(n-2)\coth r\big)\tanh r+(n-1)=1-\tanh^2r=\mathrm{sech}^2r .
\]
By \eqref{eq:gaugen}, $d\theta/dr=-\tanh r$, so $\theta=-\log\cosh r$. Along the hypersurface, $x_0+x_1=\cosh r\,e^{\theta}=1$, that is, $\la x,p\ra=-1$ for the null vector $p=e_0-e_1$. This is a horosphere centred at the ideal point $[p]$. Conversely, let $p$ be a null vector. On $\{\la x,p\ra=-1\}$ the unit normal is $\nu=p-x$, and $H\equiv1$. The boost $X_{\rm unit}=L_{01}$ satisfies $L_{01}p=-p$, so
\[
\la X_{\rm unit},\nu\ra=\la L_{01}x,p\ra=-\la x,L_{01}p\ra=\la x,p\ra=-1 .
\]
Hence $|\la X_{\rm unit},\nu\ra|\equiv1$, and $H\la\lambda X_{\rm unit},\nu\ra=\lambda$ up to orientation. This equals $1$ exactly when $\lambda=1$.

\emph{Case $\lambda>1$.} Near $r=0$, $v_\lambda<\tanh r$, since $\frac1{\lambda}<1$. Suppose the two graphs touched for the first time at some $r_1>0$. Then $v_\lambda'(r_1)\ge(\tanh)'(r_1)$, while
\[
v_\lambda'(r_1)=F_\lambda(r_1,v_\lambda(r_1))<F_{1}(r_1,\tanh r_1)=(\tanh)'(r_1),
\]
a contradiction. Hence $0<v_\lambda<\tanh r$ for all $r>0$. (At $v=0$ we have $v'>0$, so $v_\lambda$ stays positive.)

For the asymptotics we argue in two steps. First, $\sqrt{1-v^2}>\mathrm{sech}\,r$, so by \eqref{eq:driftbound} $v'\le-(n-1)v+\frac{n-1}\lambda$. Hence $\limsup v\le\frac1{\lambda}<1$. In particular, $\sqrt{1-v^2}$ is bounded away from zero, the forcing is $O(e^{-r})$, and $v=O(e^{-r})$. Second, since $\tanh r+(n-2)\coth r=n-1+O(e^{-2r})$ and $\frac{e^r}{\cosh r}=2+O(e^{-2r})$, the function $u=e^rv$ satisfies
\[
u'=-(n-2)u+\frac{2(n-1)}\lambda+O(e^{-2r})=-(n-2)(u-\kappa_n)+O(e^{-2r}).
\]
Since $n-2>0$, $u\to\kappa_n$, that is, $v=\kappa_ne^{-r}(1+o(1))$. Then $H=1/(\lambda\cosh r\cos\sigma)=\frac2\lambda e^{-r}(1+o(1))$, and
\[
\theta_\infty-\theta(r)=\int_r^\infty\frac{v}{\cosh\varrho\sqrt{1-v^2}}\,d\varrho=\int_r^\infty2\kappa_ne^{-2\varrho}(1+o(1))\,d\varrho=\kappa_ne^{-2r}(1+o(1)).
\]
The hyperplane $\{\theta=\theta_\infty\}=\Hy^n\cap\{x_1\cosh\theta_\infty-x_0\sinh\theta_\infty=0\}$ is the intersection of $\Hy^n$ with a timelike hyperplane through the origin. Hence it is totally geodesic and orthogonal to the axis.

\emph{Case $\lambda<1$.} The same comparison gives $v_\lambda>\tanh r$ for $r>0$. If $v_\lambda<1$ for all $r$, then $\sqrt{1-v^2}<\mathrm{sech}\,r$, so
\[
v'\ge-\big(\tanh r+(n-2)\coth r\big)+\frac{n-1}\lambda\longrightarrow(n-1)\Big(\frac1\lambda-1\Big)>0 .
\]
Hence $v$ exceeds $1$ at a finite $r$, a contradiction. So $v_\lambda$ reaches $1$ at some $r_+(\lambda)$, where Theorem~\ref{thm:fronts} applies. If $\lambda_1<\lambda_2$ and $r_1=r_+(\lambda_1)$, then $v_{\lambda_1}>v_{\lambda_2}$ on $(0,r_1]$ by $\partial_\lambda F_\lambda<0$, and \eqref{eq:flux} gives $\lambda_2\cc\s^{n-2}v_{\lambda_2}(r_1)<\lambda_1\cc\s^{n-2}$ at $r_1$, so $v_{\lambda_2}(r_1)<1$. Hence $r_+$ is increasing.

\emph{The limit $\lambda\uparrow1$.} Fix $R>0$. The cap $v_\lambda$ is the fixed point of a contraction that depends continuously on $\lambda$ (Lemma~\ref{lem:cap}), so $v_\lambda\to\tanh r$ uniformly on $[0,R]$ as $\lambda\to1$. Since $\tanh r\le\tanh R<1$ there, we get $v_\lambda<1$ on $[0,R]$ for $\lambda$ close to $1$, that is, $r_+(\lambda)>R$.
\qed

\subsection{Outer ends}
The following proposition contains the analysis of all outer ends in hyperbolic space, including the three-dimensional helicoidal ones.

\begin{proposition}\label{prop:outer}
Let $n\ge3$, $k,\lambda>0$, and consider the equation $v'=D(r)v+\mathcal F(r)/\sqrt{1-v^2}$ on $[R_0,\infty)$, where
\[
D=-(n-1)+O(e^{-2r}),\qquad \mathcal F=\frac{2(n-1)k}{\lambda}e^{-r}\big(1+O(e^{-2r})\big),\qquad \mathcal F'=-\mathcal F\big(1+O(e^{-2r})\big).
\]
Let $z_*=\mathcal F^2/(2(n-1)^2)$. Every solution defined on $[R,\infty)$, $R\ge R_0$, satisfies exactly one of the following:
\begin{enumerate}[label=(\Alph*)]
\item $v=\frac{2(n-1)k}{(n-2)\lambda}e^{-r}+O(e^{-2r})$;
\item $1-v=z_*(1+o(1))$.
\end{enumerate}
For $R$ large both alternatives occur, and there is exactly one solution $v_B$ of type (B). Moreover, a solution with $v(R)<v_B(R)$ is of type (A), and a solution with $v_B(R)<v(R)<1$ reaches $v=1$ at a finite $r$.
\end{proposition}

\begin{proof}
Let $z=1-v$. Then
\[
z'=-D(1-z)-\frac{\mathcal F}{\sqrt{z(2-z)}},\qquad z_*'=-2z_*\big(1+O(e^{-2r})\big).
\]
Fix $R_0$ so large that $\frac34(n-1)\le-D\le\frac54(n-1)$ and $z_*\le10^{-3}$ on $[R_0,\infty)$. Note that $\mathcal F=(n-1)\sqrt{2z_*}$.

\emph{Step 1.} If $16z_*\le z\le\frac14$, then $\mathcal F\le(n-1)\sqrt{z/8}$ and $\mathcal F/\sqrt{z(2-z)}\le(n-1)/\sqrt{14}<0.27(n-1)$, so $z'\ge(n-1)(\frac9{16}-0.27)>0.25(n-1)$. Since $z_*$ decreases, $z$ reaches $\frac14$, and then, since the forcing is $O(\mathcal F)$ for $z\in[\frac14,\frac12]$, $z$ reaches $\frac12$.

\emph{Step 2.} If $z<z_*/4$, then $\mathcal F/\sqrt{z(2-z)}>\mathcal F/\sqrt{z_*/2}=2(n-1)$, so $z'<-D-2(n-1)\le-\frac34(n-1)$. Since $(z/z_*)'=(z'+2z(1+O(e^{-2r})))/z_*<0$, the ratio $z/z_*$ decreases, the inequality persists, and $z$ reaches $0$ at a finite $r$, which is a singular orbit.

\emph{Step 3.} If $v\le\frac12$ at some large $r_1$, then $v$ stays below $\frac12$, since $v'<0$ at $v=\frac12$, and $v$ is bounded below by $\min\{v(r_1),0\}>-1$, since $v'>0$ for $v<0$. Hence the forcing is $O(e^{-r})$ and, by Gr\"onwall, $v=O(e^{-r})$. Then $u=e^rv$ satisfies $u'=(1+D)u+e^r\mathcal F(1+O(v^2))=-(n-2)\big(u-\tfrac{2(n-1)k}{(n-2)\lambda}\big)+O(e^{-2r})$, which gives (A).

\emph{Step 4.} Otherwise $v>\frac12$ for all large $r$, and by Steps 1--2, $\frac14\le\xi:=z/z_*\le16$. Since $\mathcal F/\sqrt{2z}=(n-1)/\sqrt\xi$, we get $\xi'=\frac{n-1}{z_*}(1-\xi^{-1/2})+O(1)$. As $1/z_*\to\infty$, any deviation $|\xi-1|\ge\eta$ at a large $r$ would force $\xi$ out of $[\frac14,16]$. Hence $\xi\to1$, which is (B).

\emph{Step 5 (uniqueness).} For two solutions of type (B), $z_1<z_2$, the difference $d=z_2-z_1$ satisfies $d'\ge(D+c\,\mathcal F^{-2})d$, because $(z(2-z))^{-1/2}$ has derivative of size $\gtrsim z^{-3/2}\gtrsim\mathcal F^{-3}$ on $\{z\le16z_*\}$. So $d$ grows superexponentially, while $d\le16z_*\to0$, a contradiction.

\emph{Step 6 (existence and threshold).} On $\{z(R)\in(0,\frac12)\}$, the sets of data for which $z$ reaches $0$ or $\frac12$ are open, disjoint and nonempty, by Steps 2 and 1. By connectedness, some datum is in neither, and its solution $z_B=1-v_B$ is of type (B), by Step~4. Two solutions of a scalar equation do not cross. Hence a solution with $v(R)<v_B(R)$ stays below $v_B$, so it never reaches $v=1$; it stays above $-1$ by the argument of Step~3, so it is defined on $[R,\infty)$, and by Step~5 it is not of type (B). By Steps~3 and~4 it is of type (A). A solution with $v_B(R)<v(R)<1$ stays above $v_B$. It cannot be of type (A), since then it would cross $v_B$, nor of type (B). Hence it is not defined on $[R,\infty)$, and it reaches $v=1$ at a finite $r$.
\end{proof}

\begin{proof}[Proof of Theorem~\ref{thm:outer}]
In the spherical type, $D=-(\tanh r+(n-2)\coth r)$ and $\mathcal F=(n-1)/(\lambda\cosh r)$. In the hyperbolic type, $D=-(\coth r+(n-2)\tanh r)$ and $\mathcal F=(n-1)/(\lambda\sinh r)$. In both cases the hypotheses of Proposition~\ref{prop:outer} hold with $k=1$, and $z_*=\frac2{\lambda^2}e^{-2r}(1+O(e^{-2r}))$.

By \eqref{eq:HXn} and \eqref{eq:soliton}, $H=1/(\lambda\rho\cos\sigma)=\mathcal F/((n-1)\cos\sigma)$, and $d\vartheta/dr=-\tan\sigma/\rho$, with $\rho\sim\frac12e^r$. In case (A), $H\to0$, $d\vartheta/dr=O(e^{-2r})$, so $\vartheta\to\vartheta_\infty$ with error $O(e^{-2r})$, and the distance to $\{\vartheta=\vartheta_\infty\}$ is at most $\rho|\vartheta-\vartheta_\infty|=O(e^{-r})$. The hypersurface $\{\vartheta=\vartheta_\infty\}$ is the intersection of $\Hy^n$ with a hyperplane through the origin, namely $\{x_1\cosh\theta_\infty-x_0\sinh\theta_\infty=0\}$ in the spherical type and $\{x_{n-1}\sin\varphi_\infty-x_n\cos\varphi_\infty=0\}$ in the hyperbolic type, so it is totally geodesic. In case (B), $H=\mathcal F/((n-1)\sqrt{z(2-z)})\to1$, and $d\vartheta/dr\approx-1/(\rho\sqrt{2z_*})=-(n-1)/(\rho\mathcal F)=-\lambda$. In both cases arc length is at least $r$, so the end is complete. For the translators with $\lambda=1$, the horosphere has $1-\tanh r=2e^{-2r}(1+o(1))=z_*(1+o(1))$, so by uniqueness it is the solution of type (B).
\end{proof}

\subsection{Axes of codimension two}
The following proposition describes the behaviour near an axis of codimension two. It applies to the hyperbolic type and, for $n=3$, to the helicoidal classes. The far axis in $\Sp^n$ has the same local structure, and it is treated in Section~\ref{subsec:sphere} in the parameter of \eqref{eq:systemt}, since there we follow one front through its cuspidal edges.

\begin{proposition}\label{prop:axis}
Let $\mu>0$ and consider $v'=D(r)v+\mathcal F(r)/\sqrt{1-v^2}$ for small $r>0$, where $rD(r)+1$ and $r\mathcal F(r)-\mu$ are smooth and $O(r^2)$. Let $G(v)=-v+\mu/\sqrt{1-v^2}$. There is $r_0>0$ such that the solutions that are defined at some $r_1\le r_0$ behave as follows when continued towards $r=0$. No such solution reaches $v=1$.
\begin{enumerate}[label=(\roman*)]
\item If $\mu<\frac12$, let $v_1<\frac1{\sqrt2}<v_2$ be the roots of $v\sqrt{1-v^2}=\mu$. There is exactly one solution $v_e$ on $(0,r_0]$ with $v_e\to v_1$ as $r\to0$. A solution with $v(r_1)>v_e(r_1)$ reaches $r=0$ with $v\to v_2$, and a solution with $v(r_1)<v_e(r_1)$ reaches $v=-1$ at some $r_->0$. The solutions that reach $r=0$ do so after a finite arc length.
\item If $\mu=\frac12$, there is exactly one solution $v_c$ on $(0,r_0]$ with $v_c-\frac1{\sqrt2}=O(r)$, and in fact $v_c-\frac1{\sqrt2}=O(r^2)$. A solution with $v(r_1)>v_c(r_1)$ reaches $r=0$ after a finite arc length, with
\[
v-\frac1{\sqrt2}=\frac{1+o(1)}{2\sqrt2\,|\log r|},
\]
and a solution with $v(r_1)<v_c(r_1)$ reaches $v=-1$ at some $r_->0$.
\item If $\mu>\frac12$, every solution reaches $v=-1$ at some $r_->0$. Continued as a front, in the parameter $t$ of \eqref{eq:systemt}, it reaches $r=0$ after infinitely many cuspidal points, which accumulate at $r=0$, and its total length is finite.
\end{enumerate}
\end{proposition}

\begin{proof}
With $\varsigma=-\log r$, the equation reads
\begin{equation}\label{eq:axisvs}
\frac{dv}{d\varsigma}=-G(v)-E(r,v),\qquad E=(1+rD)v+(r\mathcal F-\mu)(1-v^2)^{-1/2},
\end{equation}
so that $|E|\le Cr^2(1-v^2)^{-1/2}$, and $|E|+|\partial_vE|\le C_\delta r^2$ for $|v|\le1-\delta$. Since $G\ge\frac\mu2(1-v^2)^{-1/2}$ near $v=\pm1$, we may fix $\delta>0$ and $r_0$ such that $dv/d\varsigma\le-\frac\mu4(1-v^2)^{-1/2}$ for $1-\delta\le|v|<1$ and $r\le r_0$. Hence, as $\varsigma$ increases, no solution reaches $v=1$, and a solution that enters $\{v\le-1+\delta\}$ reaches $v=-1$ at a finite $\varsigma$. Two solutions of \eqref{eq:axisvs} do not cross, so the order between solutions is preserved.

The zeros of $G$ in $(0,1)$ satisfy $v\sqrt{1-v^2}=\mu$, and $G>0$ on $(-1,0]$. Since $\max v\sqrt{1-v^2}=\frac12$, there are no zeros when $\mu>\frac12$, a double zero at $\frac1{\sqrt2}$ when $\mu=\frac12$, with $G''(\frac1{\sqrt2})=4\sqrt2$, and two zeros $v_1<v_2$ when $\mu<\frac12$. At a zero,
\[
G'(v_i)=\frac{2v_i^2-1}{1-v_i^2},
\]
so $G'(v_1)<0<G'(v_2)$. Moreover $G>0$ on $(-1,v_1)\cup(v_2,1)$ and $G<0$ on $(v_1,v_2)$.

(i) Let $\mathcal U$ be the set of solutions that converge to $v_2$ as $\varsigma\to\infty$, and let $\mathcal L$ be the set of solutions that reach $v=-1$. We identify a solution with its value at a fixed $\varsigma_1$. Then $\mathcal L$ is open, and $\mathcal U$ is open by Lemma~\ref{lem:atr}, applied to $w=v-v_2$ with $k=G'(v_2)$. Fix $\delta'>0$ small. On $[v_1+\delta',v_2-\delta']$ we have $-G\ge c(\delta')>0$, and on $[v_2+\delta',1)$ we have $-G\le-c(\delta')$. Hence, if $\varsigma$ is large enough, depending on $\delta'$, a solution with $v(\varsigma)\ge v_1+\delta'$ enters the $\delta'$-neighbourhood of $v_2$ and belongs to $\mathcal U$. Similarly, $-G\le-c(\delta')$ on $(-1,v_1-\delta']$, so a solution with $v(\varsigma)\le v_1-\delta'$, $\varsigma$ large, belongs to $\mathcal L$. In particular $\mathcal U$ and $\mathcal L$ are nonempty, $\mathcal U$ is an upper set and $\mathcal L$ is a lower set. By connectedness there is a solution $v_e$ in neither, and by what we have just proved $|v_e-v_1|<\delta'$ for $\varsigma$ large, for every $\delta'$. So $v_e\to v_1$. If $v_a<v_b$ were two solutions converging to $v_1$, their difference $d=v_b-v_a$ would satisfy
\[
\frac{dd}{d\varsigma}=-\big(G'(\zeta)+\partial_vE\big)d\ \ge\ \tfrac12|G'(v_1)|\,d
\]
for $\varsigma$ large, where $\zeta$ lies between $v_a$ and $v_b$. Then $d$ would grow exponentially, which contradicts $d\to0$. So $v_e$ is unique, a solution above $v_e$ belongs to $\mathcal U$, and a solution below $v_e$ belongs to $\mathcal L$. The arc length $\int dr/\cos\sigma$ of the solutions that reach $r=0$ is finite, since $\cos\sigma$ has a positive limit.

(ii) Let $w=v-\frac1{\sqrt2}$. Then
\[
\frac{dw}{d\varsigma}=-2\sqrt2\,w^2+R,\qquad |R|\le C\big(|w|^3+e^{-2\varsigma}\big)\quad\text{for }|v|\le1-\delta,
\]
and there are $0<g\le g'$ and $\eta>0$ such that $G\ge g\,w^2$ on $(-1,1)$ and $G\le g'w^2$ for $|w|\le\eta$. Fix $K$ with $gK^2>2C$.

\emph{The lower region.} Let $w(\varsigma_1)\le-Ke^{-\varsigma_1}$. While $w\le-Ke^{-\varsigma}$ we have $dw/d\varsigma\le-g\,w^2+Ce^{-2\varsigma}\le-\frac g2w^2<0$, whereas $-Ke^{-\varsigma}$ increases. So the region $\{w\le-Ke^{-\varsigma}\}$ is invariant, $w\le w(\varsigma_1)<0$, and $dw/d\varsigma\le-\frac g2w(\varsigma_1)^2$. Hence $v$ reaches $-1$ at a finite $\varsigma$.

\emph{The upper region.} On the curve $w=Ke^{-\varsigma}$ we have $dw/d\varsigma\ge-(g'K^2+C)e^{-2\varsigma}>-Ke^{-\varsigma}$ for $\varsigma\ge\varsigma_0$. So the region $\{w\ge Ke^{-\varsigma}\}$ is invariant for $\varsigma\ge\varsigma_0$. In it, $dw/d\varsigma\le-\frac g2w^2$, so $w\to0$, and, once $w\le\eta$,
\[
\frac{dw}{d\varsigma}\ \ge\ -\big(2\sqrt2+C\eta+CK^{-2}\big)w^2=:-a'w^2 .
\]
Hence $1/w\le a'\varsigma+C'$, so $w\ge\alpha/\varsigma$ for $\varsigma$ large, for every $\alpha<1/a'$. Lemma~\ref{lem:sn}, with $a=2\sqrt2$, gives $\varsigma\,w\to\frac1{2\sqrt2}$. The arc length is finite, since $\cos\sigma\to\frac1{\sqrt2}$.

\emph{The threshold.} Let $\mathcal U$ and $\mathcal L$ be the sets of solutions that enter the upper and the lower region, respectively. They are open, disjoint and nonempty, $\mathcal U$ is an upper set and $\mathcal L$ is a lower set. By connectedness there is a solution $v_c$ in neither; it satisfies $|w_c|<Ke^{-\varsigma}$ for $\varsigma\ge\varsigma_0$. If $v_a<v_b$ were two such solutions, their difference $d$ would satisfy $|d\log d/d\varsigma|=|G'(\zeta)+\partial_vE|\le C''e^{-\varsigma}$, because $G'(\frac1{\sqrt2})=0$ and $|\zeta-\frac1{\sqrt2}|<Ke^{-\varsigma}$. Then $d$ would be bounded away from zero, which contradicts $d<2Ke^{-\varsigma}$. So $v_c$ is unique, the solutions above it belong to $\mathcal U$ and the solutions below it belong to $\mathcal L$. Finally, $w_c\to0$ gives
\begin{equation}\label{eq:wc}
w_c(\varsigma)=\int_\varsigma^\infty\big(G(v_c)+E(r,v_c)\big)\,d\varsigma' ,
\end{equation}
and $G(v_c)\le g'K^2e^{-2\varsigma}$, $|E|\le Ce^{-2\varsigma}$. Hence $w_c=O(e^{-2\varsigma})=O(r^2)$.

(iii) Here $G\ge g_0>0$ on $(-1,1)$. For $r$ small, $dv/d\varsigma\le-g_0/2$, so $v$ reaches $-1$ at a finite $\varsigma$. For the continuation we use the system \eqref{eq:systemt}, in which
\[
\frac{dr}{dt}=\cos^2\sigma,\qquad r\,\frac{d\sigma}{dt}=r\mathcal F+rD\sin\sigma\cos\sigma=\mu-\tfrac12\sin2\sigma+O(r^2)\ \in\ \big[c,\ \mu+1\big]
\]
for $r\le r_0$, where $c=\tfrac12\big(\mu-\tfrac12\big)>0$. As $t$ decreases, $r$ does not increase and $\sigma$ decreases. If $r$ had a positive limit $r_\infty$, then $d\sigma/dt$ would be bounded, the solution would be defined for all $t\to-\infty$, and over each half-turn of $\sigma$ the radius would decrease by $\int\cos^2\sigma\,(d\sigma/dt)^{-1}d\sigma\ge\frac{\pi r_\infty}{2(\mu+1)}$, which is absurd. So $r\to0$. Along the front,
\[
\frac{d\log r}{d\sigma}=\frac{\cos^2\sigma}{\mu-\frac12\sin2\sigma+O(r^2)}\ \in\ \Big[\frac{\cos^2\sigma}{\mu+1},\ \frac{\cos^2\sigma}{c}\Big],
\]
so over each full turn of $\sigma$ the quantity $\log r$ decreases by at least $\pi/(\mu+1)$ and at most $\pi/c$. Hence $\sigma\to-\infty$, the radii $r_k$ at successive turns decay geometrically, and the length $\int|\cos\sigma|\,dt\le c^{-1}\int r\,d\sigma$ of the $k$-th turn is $O(r_k)$. So the total length is finite. Each passage of $\sigma$ through $\frac\pi2\pmod\pi$ is a cuspidal point, by Theorem~\ref{thm:fronts}.
\end{proof}

\begin{figure}[tp]
\centering
\includegraphics[width=0.72\textwidth]{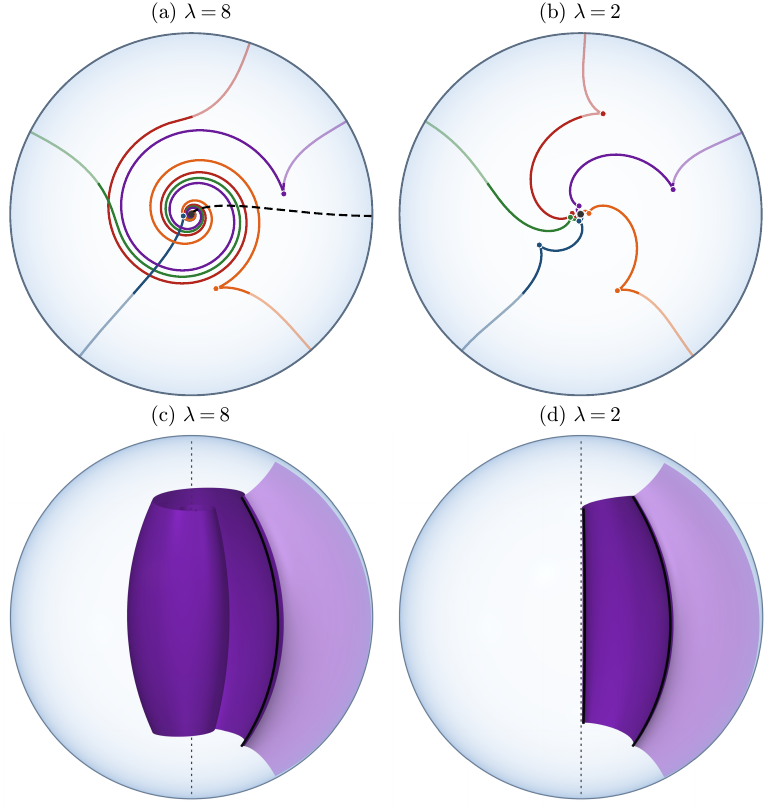}
\caption{Boost-invariant rotators of $\Hy^3$, the hyperbolic type for $n=3$ (Theorem~\ref{thm:hyp} and Corollary~\ref{cor:boost3}). (a), (b) Profiles in the totally geodesic plane orthogonal to the axis, seen as the equatorial disc of the Poincar\'e ball; the axis is the centre. Each coloured profile starts at distance $1.2$ from the axis; the inward part is dark and the outward part is light. For $\lambda=8>4$ the inward ends reach the axis along logarithmic spirals of angle $\sigma_2$, some of them after a cuspidal point; the dashed curve is the exceptional solution $v_e$, which reaches the axis along the spiral of angle $\sigma_1=\frac\pi2-\sigma_2$. For $\lambda=2<4$ every inward end reaches a singular hypercycle, and the front continues through cuspidal points that accumulate at the axis. The outward ends are complete and approach the ideal boundary. (c), (d) The surfaces generated by the purple profiles, in the Poincar\'e ball, for $|\theta|\le1.4$; the axis is dotted. In (c) the surface winds around the axis; in (d) the front is drawn between its first two inner cuspidal edges.}
\label{fig:boost}
\end{figure}

\subsection{Proof of Theorem~\ref{thm:hyp}}
Here $\rho=\s=\sinh r$, $\beta=\cosh r$, and
\begin{equation}\label{eq:hyp}
\frac{dv}{dr}=-\big(\coth r+(n-2)\tanh r\big)v+\frac{n-1}{\lambda\sinh r\sqrt{1-v^2}} .
\end{equation}

\emph{(i) Outer ends.} This is Theorem~\ref{thm:outer}. For the existence of bowl ends, let $w=\frac{M}{\lambda}e^{-r}$ with $M=\frac{10(n-1)}{n-2}$, and take $R\ge1$ with $\tanh R\ge\frac34$ and $w(R)\le\frac12$. Using $\sinh r\ge e^r/4$, $\sqrt{1-w^2}\ge\frac{\sqrt3}2$ and $\coth r+(n-2)\tanh r\ge\frac34(n-1)+\frac14$, we get, for $r\ge R$,
\[
-\big(\coth r+(n-2)\tanh r\big)w+\frac{n-1}{\lambda\sinh r\sqrt{1-w^2}}\le\Big(-\frac{3n-2}4M+\frac{8(n-1)}{\sqrt3}\Big)\frac{e^{-r}}{\lambda}<-\frac M\lambda e^{-r}=w' .
\]
So $w$ is a supersolution, while $v=0$ is a subsolution. Hence $0\le v\le w$ on $[R,\infty)$, and the solution is a bowl end.

\emph{(ii) and (iv).} Near the axis, $r\,D=-r\coth r-(n-2)r\tanh r=-1+O(r^2)$ and $r\mathcal F=\frac{n-1}\lambda\frac{r}{\sinh r}$. So Proposition~\ref{prop:axis} applies with $\mu=\frac{n-1}\lambda$, and $\mu<\frac12$ if and only if $\lambda>2(n-1)$. In case (ii), $\sin2\sigma_i=2v_i\sqrt{1-v_i^2}=2(n-1)/\lambda$. In the polar coordinates $(r,\varphi)$ of the orbit space $dr^2+\sinh^2r\,d\varphi^2$, the profile meets the radial direction at the constant limit angle $\sigma_i\neq0$, so $\varphi=\tan\sigma_i\log r+O(1)$ and the profile is a logarithmic spiral. Moreover, $H=1/(\lambda\sinh r\cos\sigma)\sim1/(\lambda r\cos\sigma_i)$. If an inward solution with $v<v_e$ reaches $v=-1$ at $r_-$, the front continues through the cuspidal point, by Theorem~\ref{thm:fronts}. By the symmetry $\sigma\mapsto\sigma+\pi$ of \eqref{eq:systemt}, the continuation is the inward solution $\tilde v=-v$ of \eqref{eq:scalar} with $\tilde v(r_-)=1>v_e(r_-)$, which reaches the axis with $\tilde v\to v_2$ by Proposition~\ref{prop:axis}(i). Case (iv) is Proposition~\ref{prop:axis}(iii); the angle $\varphi$ changes over each turn of $\sigma$ by a nonzero amount, computed as in the proof of Theorem~\ref{thm:sphere} below, so the front winds infinitely many times around the axis.

\emph{(iii) The borderline $\lambda=2(n-1)$.} Now $\mu=\frac12$, and Proposition~\ref{prop:axis}(ii) gives the solution $v_c$ and the trichotomy. Here
\[
E(r,v)=\big(1-r\coth r-(n-2)r\tanh r\big)v+\Big(\frac{r}{\sinh r}-1\Big)\frac{1}{2\sqrt{1-v^2}} ,
\]
and the expansions $r\coth r=1+\frac{r^2}3+O(r^4)$, $r\tanh r=r^2+O(r^4)$ and $\frac r{\sinh r}=1-\frac{r^2}6+O(r^4)$ give
\[
E\Big(r,\frac1{\sqrt2}+w\Big)=-\frac{2n-3}{2\sqrt2}\,r^2+O\big(r^4+r^2|w|\big).
\]
Since $w_c=O(r^2)$, we have $G(v_c)=O(r^4)$, and \eqref{eq:wc} yields
\[
w_c=-\frac{2n-3}{2\sqrt2}\int_\varsigma^\infty e^{-2\varsigma'}d\varsigma'+O(r^4)=-\frac{2n-3}{4\sqrt2}\,r^2+O(r^4).
\]
Finally, since $r\coth r>1$ and $r<\sinh r$, we have $E<0$ for $v>0$ and every $r>0$. So $dw/d\varsigma=-E>0$ at $w=0$, and $dv/d\varsigma\to-\infty$ as $v\to1$. Hence an inward solution with $v(r_1)\in(\frac1{\sqrt2},1)$ at some $r_1>0$ satisfies $v\in(\frac1{\sqrt2},1)$ for all $r<r_1$. In particular $v>v_c$ near the axis, and the solution reaches the axis.
\qed

\begin{figure}[tp]
\centering
\includegraphics[width=0.72\textwidth]{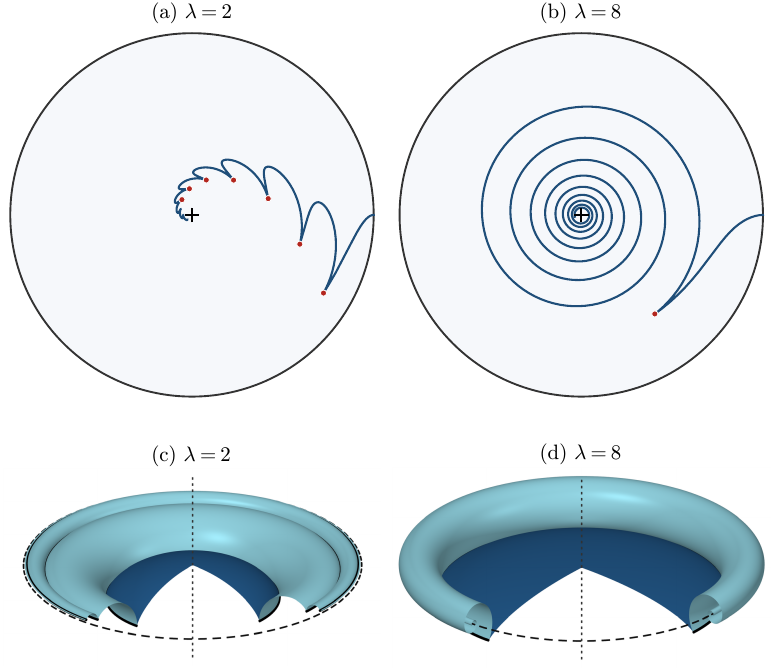}
\caption{Rotational caps of $\Sp^3$ (Theorem~\ref{thm:sphere} and Corollary~\ref{cor:sphere3}). (a), (b) Profiles in the orbit space of $\partial_\varphi$, a hemisphere whose boundary is the axis and whose centre ($+$) is the far axis, with radial coordinate $\big((\frac\pi2-r)/\frac\pi2\big)^{1/5}$. Each cap leaves the axis orthogonally, reaches its singular circle at $r_+(\lambda)$, and continues as a front. For $\lambda=2<4$ infinitely many cuspidal points accumulate at the far axis; for $\lambda=8>4$, after one cuspidal point, the front winds around the far axis along a logarithmic spiral. (c), (d) The corresponding surfaces under the stereographic projection, cut open ($0\le\varphi\le\frac{3\pi}2$): the cap (dark), the continued front (light), the axis (dotted line) and the far axis (dashed unit circle).}
\label{fig:sphere}
\end{figure}

\subsection{Proof of Theorem~\ref{thm:sphere}}\label{subsec:sphere}
Let $\varepsilon=1$, so that $D=\tan r-(n-2)\cot r$ and $\mathcal F=(n-1)/(\lambda\cos r)$. Let $v$ be the cap solution. Since $v'>0$ whenever $v=0$, we have $v>0$. For $r\ge r_n:=\arctan\sqrt{n-2}$, $\tan r\ge(n-2)\cot r$, so
\[
\frac{dv}{dr}\ge\frac{n-1}{\lambda\cos r} .
\]
If $v<1$ up to $\pi/2$, then $1>v(r_n)+\frac{n-1}\lambda\int_{r_n}^r\frac{d\varrho}{\cos\varrho}\to\infty$, which is absurd. So $v$ reaches $1$ at some $r_+<\pi/2$, which is a singular sphere by Theorem~\ref{thm:fronts}.

Beyond the singular sphere we use the parameter $t$ of \eqref{eq:systemt}, together with $x=\frac\pi2-r$:
\begin{equation}\label{eq:farfield}
\dot x=-\cos^2\sigma,\qquad \dot\sigma=\sin\sigma\cos\sigma\big(\tan r-(n-2)\cot r\big)+\frac{n-1}{\lambda\cos r}.
\end{equation}
Since $x$ is monotone and bounded, it has a limit $x_\infty$. Suppose $x_\infty>0$. Then the $\omega$-limit set in the compact cylinder $\{x=x_\infty\}\times\R/\pi\mathbb Z$ is nonempty and invariant, and it is contained in $\{\dot x=0\}=\{\cos\sigma=0\}$. But $\dot\sigma=(n-1)/(\lambda\cos r)\neq0$ there, a contradiction. Hence $x\to0$, and the front reaches the far axis.

Near $x=0$, $\tan r-(n-2)\cot r=\frac1x-\big(n-\frac53\big)x+O(x^3)$ and $\frac1{\cos r}=\frac1x+\frac x6+O(x^3)$, so
\[
\dot\sigma=\frac{B(\sigma)}x+x\,p(\sigma)+O(x^3),\qquad B(\sigma)=\frac{n-1}\lambda+\frac12\sin2\sigma,\qquad p(\sigma)=-\frac{3n-5}6\sin2\sigma+\frac{n-1}{6\lambda}.
\]

\emph{Case $\lambda<2(n-1)$.} Then $B\ge c=\frac{n-1}\lambda-\frac12>0$. Since $d\sigma/dx\le-c/x+O(x)$, $\sigma$ increases at least by $c\log(x_0/x)$. Over each full turn of $\sigma$,
\[
\log x_{k+1}-\log x_k=-\int\frac{\cos^2\sigma}{B(\sigma)}\,d\sigma\,(1+o(1))\le-\frac{\pi}{\max B}.
\]
So the radii $x_k$ at successive turns decay geometrically, and the length of the $k$-th turn is $O(x_k)$. Since $d\theta/d\sigma=-\sin\sigma\cos\sigma/B(\sigma)\,(1+O(x^2))$, the angle $\theta$ around the far axis changes per turn by the nonzero amount
\[
\Big|\int_0^{2\pi}\frac{\sin\sigma\cos\sigma}{B(\sigma)}\,d\sigma\Big|=2\pi\Big(\frac{1}{\sqrt{1-\lambda^2/(4(n-1)^2)}}-1\Big).
\]
Each passage of $\sigma$ through $\frac\pi2+k\pi$ is a cuspidal edge by Theorem~\ref{thm:fronts}.

\emph{Case $\lambda>2(n-1)$.} The zeros of $B$ in $(\frac\pi2,\pi)$ are $\sigma_*=\frac\pi2+\frac12\arcsin\frac{2(n-1)}\lambda$ and $\sigma_{**}=\pi-\frac12\arcsin\frac{2(n-1)}\lambda$, and $B$ is $\pi$-periodic. Since $B'(\sigma_*)=\cos2\sigma_*<0$, the root $\sigma_*$ is attracting. Near $\sigma_{**}$, the equation in $\varsigma=-\log x$ is $dw/d\varsigma=kw+R$ with $k>0$, and at most one solution stays near $\sigma_{**}$ for all time. Every other orbit converges to $\sigma_*\pmod\pi$ by Lemma~\ref{lem:atr}. Then $d\theta/dx=-\tan\sigma/x\,(1+O(x^2))$, so $\theta=-\tan\sigma_*\log x+O(1)\to\infty$, and the front arrives along a logarithmic spiral.

\emph{Case $\lambda=2(n-1)$.} Then $B(\sigma)=\frac12(1+\sin2\sigma)\ge0$ vanishes only at $\sigma_0=\frac{3\pi}4\pmod\pi$, and $p(\sigma_0)=\frac{2n-3}4$. Let $\varsigma=-\log x$, so that $\dot\varsigma=\cos^2\sigma/x$, and let $z=\sigma_0-\sigma$. Since $B=\sin^2z$ and $\cos^2\sigma=\frac12(1-\sin2z)$,
\[
\frac{dz}{d\varsigma}=-\frac{2\big(\sin^2z+e^{-2\varsigma}p(\sigma)+O(e^{-4\varsigma})\big)}{1-\sin2z}=-2z^2+R,\qquad |R|\le C\big(|z|^3+e^{-2\varsigma}\big).
\]
This is the same equation as for $n=3$, up to the value of $p(\sigma_0)$, which only enters the constant $C$. Since $d\sigma/d\varsigma>0$ away from $\sigma_0$, each approach of $\sigma$ to $\sigma_0\pmod\pi$ enters $\{0<z\le\eta\}$ at some $\varsigma_e$. While $z\ge0$, the inequality $z'\ge-2(1+C\eta)z^2-Ce^{-2\varsigma}$ and comparison with the unforced Riccati equation give, at $\varsigma_1=2\varsigma_e$,
\[
z(\varsigma_1)\ge\frac{1}{2(1+C\eta)(\varsigma_1-\varsigma_e)+1/\eta}-\frac C2e^{-2\varsigma_e}\ge\frac{\alpha}{\varsigma_1}
\]
for any fixed $\alpha<\frac1{2(1+C\eta)}$, provided $\varsigma_e$ is large. Then Lemma~\ref{lem:sn} applies, with $a=2$: the front converges to $\sigma_0$ and $z|\log x|\to\frac12$. If $\varsigma_e$ is not large, the front may cross $\sigma_0$. It then performs a further half-turn of $\sigma$, which takes $\varsigma$-time at least $\int\cos^2\sigma\,d\sigma/\max B>0$. Hence after finitely many half-turns, and so finitely many cuspidal edges, the entry time $\varsigma_e$ is large and the front is captured.
\qed

\begin{remark}
Solutions of \eqref{eq:scalar} that do not meet the axis behave differently. When $\lambda>2(n-1)$, a solution with $v<0$ near the far axis can converge to the root $v_-<-\frac1{\sqrt2}$ of $-v\sqrt{1-v^2}=\frac{n-1}\lambda$. Indeed, near the far axis, $dv/d\varsigma=v+\frac{n-1}{\lambda\sqrt{1-v^2}}+O(e^{-2\varsigma})$ with $\varsigma=-\log x$, and $v_-$ is an attracting zero of the right-hand side, as in Proposition~\ref{prop:axis}. Such a solution reaches the far axis along a spiral without ever meeting $v=1$. The statement of Theorem~\ref{thm:sphere} is therefore restricted to caps.
\end{remark}

\subsection{Proof of Theorem~\ref{thm:rplus}}

\emph{(i) Small speed.} Let $\varrho=r/\lambda$ and $\mu=\varepsilon\lambda^2$. The functions $u\,\cc(u)/\s(u)$, $\varepsilon u\,\s(u)/\cc(u)$ and $1/\cc(u)$ are analytic functions of $\varepsilon u^2$, and the second one vanishes at $u=0$. Hence, along a cap, \eqref{eq:scalar} becomes
\[
\frac{dv}{d\varrho}=\Big(-\frac{n-2}\varrho+\mu\varrho\,\mathcal A(\mu\varrho^2)\Big)v+\frac{(n-1)\mathcal B(\mu\varrho^2)}{\sqrt{1-v^2}},\qquad \mathcal B(0)=1,
\]
with $\mathcal A,\mathcal B$ analytic. For $\mu=0$ this is the equation of the rotational IMCF translators of $\R^n$ with speed $1$. Indeed, for a hypersurface of revolution with profile angle $\sigma$ and $v=\sin\sigma$, $(n-1)H=v'+(n-2)v/\varrho$ and $\la e,\nu\ra=\sqrt{1-v^2}$, where $e$ is the direction of the axis.

The contraction of Lemma~\ref{lem:cap} depends smoothly on $\mu$, and so does the cap. In the parameter $t$ of \eqref{eq:systemt}, the cap solves a smooth vector field depending smoothly on $\mu$, and it crosses $\{\sigma=\frac\pi2\}$ transversally, with $d\sigma/dt=\mathcal F_0>0$. By the implicit function theorem, $\varrho_+(\mu)=r_+/\lambda$ is a smooth function of $\mu$. Hence $r_+=\lambda\varrho_+(\varepsilon\lambda^2)=c_0(n)\lambda+O(\lambda^3)$, with $c_0(n)=\varrho_+(0)$.

To bound $c_0(n)$, note that for $\mu=0$ the flux identity reads $\varrho^{n-2}v=(n-1)\int_0^\varrho t^{n-2}(1-v^2)^{-1/2}dt\ge\varrho^{n-1}$. The inequality is strict for $\varrho>0$, so $v>\varrho$ and $c_0(n)<1$. Conversely, while $v\le\frac1{\sqrt2}$ we have $\varrho^{n-2}v\le\sqrt2\varrho^{n-1}$, that is, $v\le\sqrt2\varrho$. By continuity, $v\le\frac1{\sqrt2}<1$ on $[0,\frac12]$, so $c_0(n)>\frac12$.

\emph{(ii) Large speed, $\varepsilon=1$.} Let
\[
I(r)=\int_0^r\sin^{n-2}\varrho\,d\varrho,\qquad W_n=I\big(\tfrac\pi2\big)=\int_0^{\pi/2}\cos^{n-2}t\,dt,\qquad V(r)=\frac{I(r)}{\cos r\,\sin^{n-2}r}.
\]
The function $V$ is increasing. Indeed, $V'$ has the sign of $J=\cos r\sin^{n-1}r+I(r)\big(\sin^2r-(n-2)\cos^2r\big)$, and $J(0)=0$, $J'=\cos^2r\sin^{n-2}r+2(n-1)I\sin r\cos r>0$. Moreover, $V=\frac{W_n+O(x)}{x}$ as $x=\frac\pi2-r\to0$. Write $\Phi=\lambda\cos r\sin^{n-2}r\,v$. By \eqref{eq:flux},
\[
\Phi(r)=(n-1)\int_0^r\frac{\sin^{n-2}\varrho}{\sqrt{1-v^2}}\,d\varrho\ \ge\ (n-1)I(r),\qquad\text{so}\qquad v\ge\frac{(n-1)V}\lambda .
\]

\emph{Step 1.} While $v\le\frac1{\sqrt2}$, we have $\Phi\le\sqrt2(n-1)I$, that is, $v\le\sqrt2(n-1)V/\lambda$. Let $r_\lambda$ be defined by $V(r_\lambda)=\lambda/(2(n-1))$. By continuity, $v\le\sqrt2(n-1)V/\lambda\le\frac1{\sqrt2}$ on $[0,r_\lambda]$. Hence $r_+>r_\lambda$, and $\frac\pi2-r_\lambda=\frac{2(n-1)W_n+O(1/\lambda)}{\lambda}$.

\emph{Step 2.} On $[r_\lambda,r_+]$ we have $v\ge(n-1)V/\lambda\ge\frac12$ and $\tan r-(n-2)\cot r\ge\frac{1}{2x}\ge c\lambda$ for $\lambda$ large. By \eqref{eq:master}, $\sigma'\ge(\tan r-(n-2)\cot r)\sin\sigma\ge\frac{c\lambda}2$. Since $\sigma\in(0,\frac\pi2)$, the arc length of this piece is $O(1/\lambda)$.

\emph{Step 3.} Write $\Phi(r_+)=(n-1)\big(I(r_+)+\delta\big)$, with $\delta=\int_0^{r_+}\sin^{n-2}\varrho\big((1-v^2)^{-1/2}-1\big)d\varrho$. On $[0,r_\lambda]$, $v\le\frac1{\sqrt2}$, so the integrand is at most $v^2\le2(n-1)^2V^2/\lambda^2$, and $\int_0^{r_\lambda}V^2=O(1/(\frac\pi2-r_\lambda))=O(\lambda)$. On $[r_\lambda,r_+]$, the integral is at most the arc length, which is $O(1/\lambda)$ by Step~2. Hence $\delta=O(1/\lambda)$.

\emph{Step 4 (first order).} At $r_+$ we have $v=1$, so, with $x_+=\frac\pi2-r_+\le\frac\pi2-r_\lambda=O(1/\lambda)$,
\[
\lambda x_+\big(1+O(x_+^2)\big)=\lambda\cos r_+\sin^{n-2}r_+=\Phi(r_+)=(n-1)W_n+O(x_+)+O(1/\lambda).
\]
Hence $\lambda x_+=(n-1)W_n+O(1/\lambda)$, and $(n-1)W_n=(n-1)\frac{\sqrt\pi\,\Gamma(\frac{n-1}2)}{2\,\Gamma(\frac n2)}=\frac{\sqrt\pi\,\Gamma(\frac{n+1}2)}{\Gamma(\frac n2)}=K_n$. In particular, by \eqref{eq:volcap}, $\mathrm{Vol}(\Sigma_\lambda)\to\omega_{n-2}W_n=\frac12\mathrm{Vol}(\Sp^{n-1})$, which is the last assertion of (iii).

\emph{Step 5 (a second identity).} We now write everything in terms of $x=\frac\pi2-r$, so that $\cos r=\sin x$, $\sin r=\cos x$, $x'=-\cos\sigma$, and \eqref{eq:master} reads
\[
\sigma'=\Big(\cot x-(n-2)\tan x\Big)\sin\sigma+\frac{n-1}{\lambda\sin x\cos\sigma}.
\]
Let $\delta(s)=\Phi/(n-1)-I(r)=\int_0^s\cos^{n-2}x\,(1-\cos\sigma)\,du$, so that the number $\delta$ of Step~3 is $\delta(s_+)$, and let $A=\sin x\,(1-\cos\sigma)$. A direct computation gives
\[
A'=\cos x\,(1-\cos\sigma)-(n-2)\frac{\sin^2x}{\cos x}\sin^2\sigma+\frac{n-1}\lambda\tan\sigma ,
\]
and therefore
\begin{equation}\label{eq:clairaut}
(\delta-A)'=\big(\cos^{n-2}x-\cos x\big)(1-\cos\sigma)+(n-2)\frac{\sin^2x}{\cos x}\sin^2\sigma-\frac{n-1}\lambda\tan\sigma .
\end{equation}
Both $\delta$ and $A$ vanish at the vertex, and at the singular sphere $A=\sin x_+$ and $\delta=\frac\lambda{n-1}\sin x_+\cos^{n-2}x_+-I(r_+)$. Since $I(r_+)=W_n-x_++O(x_+^3)$ and $x_+=O(1/\lambda)$, integrating \eqref{eq:clairaut} along the cap gives
\begin{equation}\label{eq:second}
\lambda x_+-K_n=(n-1)(\mathcal E_1+\mathcal E_2)-\frac{(n-1)^2}{\lambda}\int_0^{s_+}\tan\sigma\,ds+O(\lambda^{-2}),
\end{equation}
where $\mathcal E_1$ and $\mathcal E_2$ are the integrals of the first two terms on the right-hand side of \eqref{eq:clairaut}. In the flat model of the far axis, $\sin x\sin\sigma$ is the Clairaut constant of a straight line, and \eqref{eq:second} expresses that the length of the segment from a far point to the foot of the perpendicular through the origin differs from the radial length by a second-order quantity.

\emph{Step 6 (the error terms).} We claim that $\mathcal E_1,\mathcal E_2=O(\lambda^{-2})$. On $[0,r_\lambda]$ we have $v\le\sqrt2(n-1)V/\lambda$, $1-\cos\sigma\le\sin^2\sigma=v^2$ and $ds\le\sqrt2\,dr$. Moreover $|\cos^{n-2}x-\cos x|\le(n-3)\sin^2x\cos x$, and $\sin x\,V=I(r)/\sin^{n-2}r\le r$, because the integrand of $I$ is increasing. Hence both integrands are bounded by $C\lambda^{-2}\sin^2x\,V^2/\cos x\le C\lambda^{-2}r^2/\sin r\le C'\lambda^{-2}$. On $[r_\lambda,r_+]$ the integrands are bounded by $C\sin^2x\le C'\lambda^{-2}$, and the arc length is $O(1/\lambda)$ by Step~2.

\emph{Step 7 (the logarithm).} We claim that
\begin{equation}\label{eq:logterm}
\int_0^{s_+}\tan\sigma\,ds=\frac{2K_n\log\lambda}{\lambda}+O\Big(\frac1\lambda\Big).
\end{equation}
Together with \eqref{eq:second} and Step~6 this gives $\lambda x_+=K_n-2(n-1)^2K_n\lambda^{-2}\log\lambda+O(\lambda^{-2})$, which is (ii).

\emph{Outer part.} On $[0,r_\lambda]$, $\tan\sigma\,ds=v(1-v^2)^{-1}dr$ and $v^2\le\frac12$, so $v\le v(1-v^2)^{-1}\le v+2v^3$. By the flux identity $v=\frac{(n-1)V}\lambda\big(1+\delta/I\big)$, and $0\le\delta(r)\le I(r)\max_{[0,r]}\big((1-v^2)^{-1/2}-1\big)\le2(n-1)^2V(r)^2I(r)/\lambda^2$, since $V$ is increasing. Hence $v=\frac{(n-1)V}\lambda+O(V^3/\lambda^3)$. Now $V$ is bounded on $[0,\frac\pi4]$ and $V=\frac{W_n}x-1+O(x)$ as $x\to0$, while $\frac\pi2-r_\lambda=2K_n\lambda^{-1}(1+O(1/\lambda))$. Therefore $\int_0^{r_\lambda}V\,dr=W_n\log\lambda+O(1)$ and $\int_0^{r_\lambda}V^3dr=O(\lambda^2)$, and
\[
\int_0^{s(r_\lambda)}\tan\sigma\,ds=\frac{(n-1)W_n\log\lambda}\lambda+O\Big(\frac1\lambda\Big)=\frac{K_n\log\lambda}\lambda+O\Big(\frac1\lambda\Big).
\]

\emph{Inner part.} On $[r_\lambda,r_+]$ we use $\sigma$ as the parameter; it increases from $\sigma_\lambda\in[\frac\pi6,\frac\pi4]$ to $\frac\pi2$, since $\frac12\le v(r_\lambda)\le\frac1{\sqrt2}$ and $\sigma'>0$ by Step~2. There, $I(r)=W_n+O(1/\lambda)$ and $0\le\delta\le\delta(r_+)=O(1/\lambda)$ by Step~3, so $\Phi=K_n(1+O(1/\lambda))$, that is,
\[
\lambda\,x\,\sin\sigma=K_n\big(1+O(1/\lambda)\big)
\]
uniformly on this piece. Since $\cot x-(n-2)\tan x=x^{-1}(1+O(x^2))$ and $\sin x=x(1+O(x^2))$, the equation for $\sigma'$ becomes
\[
\sigma'=\frac{1+O(x^2)}{x}\Big(\sin\sigma+\frac{n-1}{\lambda\cos\sigma}\Big)=\frac{\sin\sigma}{K_n}\Big(\lambda\sin\sigma+\frac{n-1}{\cos\sigma}\Big)\big(1+O(1/\lambda)\big),
\]
and therefore
\[
\int_{s(r_\lambda)}^{s_+}\tan\sigma\,ds=\int_{\sigma_\lambda}^{\pi/2}\frac{\tan\sigma}{\sigma'}\,d\sigma=K_n\big(1+O(1/\lambda)\big)\int_{\sigma_\lambda}^{\pi/2}\frac{d\sigma}{\lambda\sin\sigma\cos\sigma+n-1}.
\]
With $w=\frac\pi2-\sigma\in[0,\frac\pi3]$, the last integral is $\int_0^{\pi/2-\sigma_\lambda}\big(\frac\lambda2\sin2w+n-1\big)^{-1}dw$. Replacing $\frac12\sin2w$ by $w$ changes the integrand by $O(\lambda w^3)/(c\lambda w+n-1)^2=O(w/\lambda)$, and
\[
\int_0^{\pi/2-\sigma_\lambda}\frac{dw}{\lambda w+n-1}=\frac1\lambda\log\frac{\lambda(\frac\pi2-\sigma_\lambda)+n-1}{n-1}=\frac{\log\lambda}\lambda+O\Big(\frac1\lambda\Big).
\]
So the inner part is also $K_n\lambda^{-1}\log\lambda+O(1/\lambda)$, and \eqref{eq:logterm} follows. The two halves of the logarithm come from the two sides of the scale $x\sim1/\lambda$: the outer region, where the cap is a small graph over the totally geodesic hemisphere, and the boundary layer near the singular sphere, where $\cos\sigma\sim1/\lambda$.
\qed

\section{The parabolic type}\label{sec:parab}

The map
\[
\Psi(z,h)=\Big(\frac{A+1}{2h},\frac{A-1}{2h},\frac{z}h\Big),\qquad z\in\R^{n-1},\quad A=h^2+|z|^2,
\]
is an isometry from the upper half-space $\U^n$ onto $\Hy^n$. Indeed, $x_0+x_1=A/h$ and $x_0-x_1=1/h$, which give $\la\Psi,\Psi\ra=-1$ and $\Psi^*\la\cdot,\cdot\ra=h^{-2}(|dz|^2+dh^2)$. Moreover, $\Psi_*\partial_{z_j}=B_j\Psi$, where $B_jx=x_{j+1}(e_0+e_1)+(x_0-x_1)e_{j+1}$ is a null rotation fixing $e_0+e_1$, with $B_j^3=0$ and $[B_j,B_k]=0$. So the translations in $z$ form the abelian group of parabolic isometries fixing the ideal point $[e_0+e_1]$, which corresponds to $h=\infty$.

\begin{figure}[tp]
\centering
\includegraphics[width=0.80\textwidth]{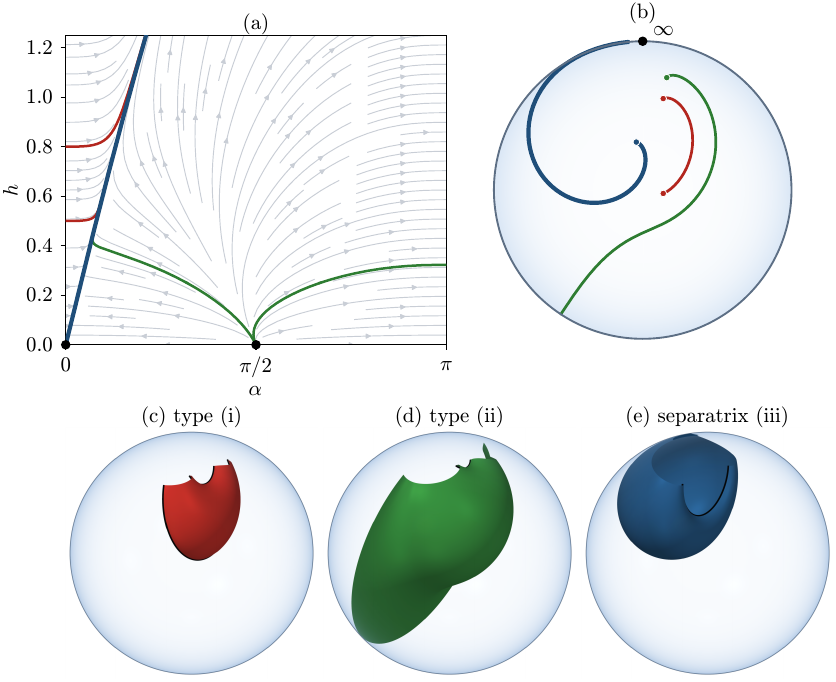}
\caption{Parabolic solitons of $\Hy^3$ (Theorem~\ref{thm:parab} and Corollary~\ref{cor:parab3}), $\lambda=2$. (a) Phase plane $(\alpha,h)$ of \eqref{eq:parab}. Red: type (i), with $\alpha\to0$ in backward time. Green: type (ii), with $(h,\alpha)\to(0,\frac\pi2)$. Thick blue: the separatrix (iii), tangent to $\alpha=\frac h2$ at the origin. (b) One profile of each type, with the same colours, in the totally geodesic plane $\{z_2=0\}$, seen as a disc of the Poincar\'e ball; the ideal point $\infty$ fixed by the parabolic group is at the top. Each profile is shown up to an isometry of $\Hy^3$ fixing $\infty$. The profile of type (ii) reaches the ideal boundary orthogonally, and the separatrix approaches $\infty$ tangentially to the ideal boundary. (c)--(e) The corresponding surfaces in the Poincar\'e ball, swept by horizontal translations in $z_2$; $\infty$ is the top point.}
\label{fig:parab}
\end{figure}

\begin{proof}[Proof of Theorem~\ref{thm:parab}]
In the orthonormal frame $E_j=h\partial_{z_j}$, $E_h=h\partial_h$, the conformal change formula gives
\[
\nabla_{E_j}E_j=E_h,\qquad \nabla_{E_j}E_h=-E_j\qquad(1\le j\le n-1),
\]
and the other covariant derivatives of frame fields vanish. The profile $T=\cos\alpha E_1+\sin\alpha E_h$ has geodesic curvature $\alpha'+\cos\alpha$ with respect to $N=\sin\alpha E_1-\cos\alpha E_h$. Each orbit direction $E_j$, $2\le j\le n-1$, is principal, with principal curvature $-\la\nabla_{E_j}E_j,N\ra=\cos\alpha$. Hence
\[
(n-1)H=\alpha'+(n-1)\cos\alpha,\qquad \la\lambda\partial_{z_1},N\ra=\frac{\lambda\sin\alpha}{h}.
\]
So \eqref{eq:soliton} reads $\alpha'=-(n-1)\cos\alpha+(n-1)h/(\lambda\sin\alpha)$. The dilation $(z,h)\mapsto\varrho(z,h)$ is an isometry that normalises the group of horizontal translations and replaces $\lambda$ by $\lambda\varrho$, so we may take $\lambda=n-1$, which gives \eqref{eq:parab}. As checks, a horosphere $h\equiv h_0$ has $\alpha\equiv0$ and $H=1$, and a vertical hyperplane has $\alpha\equiv\frac\pi2$ and $H=0$.

With $d\tau=ds/\sin\alpha$, the equations become
\[
\dot h=h\sin^2\alpha,\qquad \dot\alpha=h-(n-1)\sin\alpha\cos\alpha .
\]
The line $\{h=0\}$ is invariant. On it, $\alpha=0$ and $\alpha=\pi$ are attracting and $\alpha=\frac\pi2$ is repelling. In the interior, $\dot h>0$, so there are no equilibria.

\emph{Forward end.} Since $\dot\alpha=h>0$ on $\{\alpha=0\}$, the region $\{0<\alpha<\pi\}$ can only be left through $\alpha=\pi$. Suppose that $\alpha<\pi$ for all $\tau$. The function $h$ is increasing. If it were bounded, the $\omega$-limit set of the orbit would be a nonempty invariant subset of $\{h=h_\infty\}\cap\{\dot h=0\}=\{h=h_\infty,\ \sin\alpha=0\}$, which is impossible because $\dot\alpha=h_\infty\neq0$ there. So $h\to\infty$, and then $\dot\alpha\ge h-\frac{n-1}2\to\infty$, a contradiction. Hence $\alpha$ reaches $\pi$ at a finite $\tau$, with $h$ bounded. Near $\alpha=\pi$, the analysis of Theorem~\ref{thm:fronts}, with $\sin\alpha$ in place of $\cos\sigma$, gives a cuspidal edge along the horospherical plane $\{h=h_{\max},z_1=z_1^*\}$.

\emph{Backward end.} The point $(0,\frac\pi2)$ is an unstable node, with linearisation $\left(\begin{smallmatrix}1&0\\1&n-1\end{smallmatrix}\right)$. The point $(0,0)$ has eigenvalues $0$ and $-(n-1)$, with centre direction $\alpha=h/(n-1)$, along which $\dot h=h^3/(n-1)^2+O(h^4)>0$. By centre-manifold theory, exactly one orbit in $\{h>0\}$ has $\alpha$-limit $(0,0)$. Any other backward orbit either reaches $\alpha=0$ with $h>0$ (since $\dot\alpha=h>0$ there), which is case (i), or has $h\to0$. In the latter case its $\alpha$-limit is an invariant subset of $\{h=0\}$ different from $(0,0)$, hence it is the node, which is case (ii). Along case (ii), $h\sim Ce^s$ and $z_1'=h\cos\alpha\to0$ exponentially, so the end is complete and $H=h/((n-1)\sin\alpha)\to0$. Along case (iii), $h'\approx h^2/(n-1)$ gives $h\sim(n-1)/|s|$ as $s\to-\infty$. Hence the end is complete, $z_1\to-\infty$, and $H=h/((n-1)\sin\alpha)\to1$. The sets of data with fates (i) and (ii) are open, and they are separated by the orbit (iii).
\end{proof}

\section{The three-dimensional case}\label{sec:three}

\subsection{Commuting pairs}
\begin{proof}[Proof of Proposition~\ref{prop:pairs}]
Let $\varepsilon=1$. The Lie algebra $\mathfrak{so}(4)$ is compact of rank $2$. Every abelian subalgebra is contained in a maximal one, and the maximal abelian subalgebras are the Cartan subalgebras, which are all conjugate. Hence $\mathrm{span}\{X,Y\}$ is conjugate to $\mathfrak t=\mathrm{span}\{R_{01},R_{23}\}$, whose elements generate $\partial_\theta$ and $\partial_\varphi$.

Let $\varepsilon=-1$. We identify $x\in\R^{1,3}$ with the Hermitian matrix $\hat x=x_0\mathrm{Id}+\sum x_i\sigma_i$, where the $\sigma_i$ are the Pauli matrices. Then $\det\hat x=x_0^2-|\vec x|^2$, and the action $\hat x\mapsto S\hat xS^*$, $S\in SL(2,\mathbb C)$, gives $\mathfrak{so}(1,3)\cong\mathfrak{sl}(2,\mathbb C)$. Under this isomorphism, $\frac12\mathrm{diag}(1,-1)$ corresponds to the boost $L_{03}$, $\frac i2\mathrm{diag}(1,-1)$ to the rotation $R_{12}$, and $\mathfrak n=\left(\begin{smallmatrix}0&1\\0&0\end{smallmatrix}\right)$, $i\mathfrak n$ to the null rotations $P_1,P_2$ fixing the null direction $e_0+e_3$. Every nonzero $Z\in\mathfrak{sl}(2,\mathbb C)$ is conjugate either to $\mu\,\mathrm{diag}(1,-1)/2$, $\mu\in\mathbb C^*$, unique up to $\mu\mapsto-\mu$, or to $\mathfrak n$.

In the first case, the centraliser of $Z$ is $\mathbb C\,\mathrm{diag}(1,-1)$, that is, the real span of $L_{03}$ and $R_{12}$. For $W=\left(\begin{smallmatrix}\alpha&\beta\\\gamma&-\alpha\end{smallmatrix}\right)$ we have $[W,\mathfrak n]=\left(\begin{smallmatrix}-\gamma&2\alpha\\0&\gamma\end{smallmatrix}\right)$, so the centraliser of $\mathfrak n$ is $\mathbb C\mathfrak n=\mathrm{span}\{P_1,P_2\}$. A two-dimensional abelian subalgebra lies in the centraliser of each of its nonzero elements. Hence it equals a Cartan subalgebra if it contains a semisimple element, and it is conjugate to $\mathrm{span}\{P_1,P_2\}$ otherwise.

Conjugating the geodesic fixed by the Cartan subalgebra to $\Hy^3\cap\{x_2=x_3=0\}$ gives case (i). In the parabolic case, we send the fixed null direction to the ideal point $h=\infty$ of the half-space model (Section~\ref{sec:parab}). Then $P_1,P_2$ become $\partial_{z_1},\partial_{z_2}$. Rotations of the $(z_1,z_2)$-plane normalise $\mathrm{span}\{P_1,P_2\}$, so we may take $Y=\partial_{z_2}$.
\end{proof}

\subsection{The reduction with an arbitrary pitch}
The fields $\partial_\theta$ and $\partial_\varphi$ are Killing and commute, with $|\partial_\theta|^2=\cc^2$ and $|\partial_\varphi|^2=\s^2$. The functions $r$ and $\psi=b\theta-a\varphi$ are $Y$-invariant. The horizontal lift of $\partial_\psi$ is $\frac{b\s^2}{|Y|^2}\partial_\theta-\frac{a\cc^2}{|Y|^2}\partial_\varphi$, and its squared norm is $\s^2\cc^2/|Y|^2$. Hence, locally, the orbit space carries the metric
\begin{equation}\label{eq:quot}
g_{\rm quot}=dr^2+\rho(r)^2d\psi^2,\qquad \rho=\frac{\s\cc}{|Y|},\qquad\text{so that}\qquad \rho\,|Y|=\s\,\cc .
\end{equation}
If $\varepsilon=1$ and $a/b\notin\mathbb Q$, the orbits of $Y$ are dense in the flat tori $\{r=\mathrm{const}\}$, and the orbit space is not Hausdorff. All our arguments are local along the profile curve and remain valid in this case.

The vector field $W=-b\s^2\partial_\theta+a\cc^2\partial_\varphi$ is orthogonal to $Y$, with $|W|^2=\s^2\cc^2|Y|^2$ and $\la X,W\ra=-(pb-qa)\s^2\cc^2=-\lambda\s^2\cc^2$.

\begin{proof}[Proof of Theorem~\ref{thm:master3}]
We write the surface as $F(s,u)=x\big(r(s),\Theta(s)+au,\Phi(s)+bu\big)$, so that $F_u=Y$. We impose $\la F_s,F_u\ra=0$ and $|F_s|=1$. The first condition gives $\Theta'=-b\s^2\mu$ and $\Phi'=a\cc^2\mu$ for some function $\mu$. The second becomes $r'^2+\mu^2\s^2\cc^2|Y|^2=1$, which we solve with an angle function $\sigma$:
\begin{equation}\label{eq:gauge}
r'=\cos\sigma,\qquad \Theta'=-\frac{b\,\s\sin\sigma}{\cc|Y|},\qquad \Phi'=\frac{a\,\cc\sin\sigma}{\s|Y|}.
\end{equation}
Then $F_s=\cos\sigma\,\partial_r+\sin\sigma\,\widehat W$, where $\widehat W=W/|W|$, and $\{F_s,\widehat Y\}$ is an orthonormal frame, with $\widehat Y=Y/|Y|$. Let $N=\sin\sigma\,\partial_r-\cos\sigma\,\widehat W$.

\emph{Orbit direction.} Since $Y$ is Killing, $\overline\nabla_YY=-\frac12\mathrm{grad}|Y|^2$, so $\overline\nabla_{\widehat Y}\widehat Y=-\mathrm{grad}\log|Y|$, and $h(\widehat Y,\widehat Y)=(\log|Y|)'\sin\sigma$.

\emph{Profile direction.} The projection onto the orbit space is a Riemannian submersion, and $F_s$ and $N$ are horizontal. By O'Neill's formula and the computation in the proof of Theorem~\ref{thm:master}, with the metric \eqref{eq:quot}, $h(F_s,F_s)=\sigma'+\frac{\rho'}{\rho}\sin\sigma$.

Adding the two entries and using \eqref{eq:quot} and $(\log(\s\cc))'=\frac\cc\s-\varepsilon\frac\s\cc$, we obtain $2H=\sigma'+(\frac{\cc}{\s}-\varepsilon\frac{\s}{\cc})\sin\sigma$. Finally, since $\la X,\partial_r\ra=0$, $\la X,N\ra=-\cos\sigma\,\la X,W\ra/|W|=\lambda\,\s\cc\cos\sigma/|Y|$.
\end{proof}

For $(a,b)=(0,1)$, $\Theta'=-\sin\sigma/\cc$, which is \eqref{eq:gaugen} for the spherical type. For $(a,b)=(1,0)$, $\Phi'=\sin\sigma/\s$ and $\lambda=-q$, which is \eqref{eq:gaugen} for the hyperbolic type after $\varphi\mapsto-\varphi$. This shows that \eqref{eq:master3} agrees with \eqref{eq:master} for $n=3$.

\subsection{Proofs of the corollaries}
Corollaries~\ref{cor:transl3}, \ref{cor:sphere3}, \ref{cor:rplus3} and~\ref{cor:parab3} are Theorems~\ref{thm:transl}, \ref{thm:sphere}, \ref{thm:rplus} and~\ref{thm:parab} with $n=3$, since then $\kappa_3=4/\lambda$, $2(n-1)=4$, $K_3=2$ and $\lambda=n-1=2$. The boost-invariant rotators are the hyperbolic type for $n=3$, with $G=\{\exp(t\partial_\theta)\}$ the boosts along the axis, so Corollary~\ref{cor:boost3} is Theorem~\ref{thm:hyp} with $n=3$.

For Corollary~\ref{cor:outer3}, the equation \eqref{eq:master3} has $D=-(\tanh r+\coth r)=-2+O(e^{-4r})$ and $\mathcal F=2|Y|/(\lambda\s\cc)$. Since $|Y|=\frac{e^r}2\sqrt{a^2+b^2}(1+O(e^{-2r}))$ and $(\log|Y|)'=1+O(e^{-2r})$, the hypotheses of Proposition~\ref{prop:outer} hold with $n=3$ and $k=\sqrt{a^2+b^2}$, and $z_*=\mathcal F^2/8=\frac{\kappa^2}8e^{-2r}(1+O(e^{-2r}))$. As in the proof of Theorem~\ref{thm:outer}, $H=\mathcal F/(2\cos\sigma)$ and $d\psi/dr=-\tan\sigma/\rho$, with $\rho\mathcal F=2/\lambda$. The level set $\{\psi=\psi_\infty\}$ has $\sigma\equiv0$, hence $H\equiv0$ by \eqref{eq:HX}.

\begin{figure}[tp]
\centering
\includegraphics[width=0.80\textwidth]{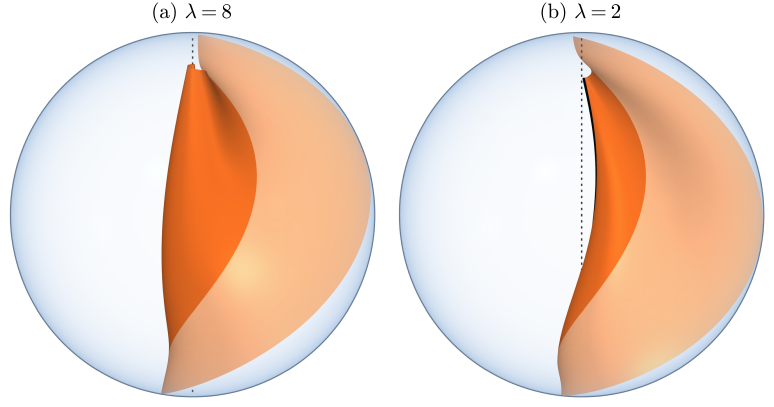}
\caption{Helicoidal solitons of $\Hy^3$ with $(a,b)=(1,1)$ (Theorem~\ref{thm:helic}), in the Poincar\'e ball. The profile starts at distance $0.8$ from the axis (dotted) with $\sin\sigma=0.6$, and the surface is swept by the screw motion generated by $Y$; the inward part is dark and the outward part is light. (a) $\lambda=8>4|a|$: the inner end reaches the axis along a logarithmic spiral. (b) $\lambda=2<4|a|$: the inner end reaches a singular orbit of $Y$ (black). In both cases the outer end is a bowl end, asymptotic to the minimal helicoid of pitch $(1,1)$.}
\label{fig:helic}
\end{figure}

\subsection{Helicoidal and Hopf classes}
\begin{proof}[Proof of Theorem~\ref{thm:helic}]
(i) The ends with $r\to\infty$ are described by Corollary~\ref{cor:outer3}. The minimal surface $\{\psi=\psi_\infty\}$ is the orbit under $Y$ of a geodesic ray orthogonal to the axis, that is, the minimal helicoid of pitch $(a,b)$ \cite{dCD}. Near $r=0$, $|Y|=|a|(1+O(r^2))$, so $r\mathcal F=\frac{2|a|}\lambda(1+O(r^2))$ and $rD=-1+O(r^2)$. Proposition~\ref{prop:axis} applies with $\mu=2|a|/\lambda$, which is the case $n=3$ of Theorem~\ref{thm:hyp} with $\lambda/|a|$ in place of $\lambda$.

(ii) When $\varepsilon=1$, the change $r\mapsto\frac\pi2-r$ exchanges $a$ and $b$, so the same holds at the far axis with $|b|$. For $r\ge\pi/4$ and $v\ge0$, we have $|Y|\ge|b|/\sqrt2$, the drift is nonnegative, and $dv/dr\ge\sqrt2|b|/(\lambda\cos r)$. As in Theorem~\ref{thm:sphere}, $v$ reaches $1$ before $\pi/2$. If $v<0$, the solution either crosses $v=0$ or is captured by the far-axis spiral, which is possible only if $\lambda\ge4|b|$. The same holds at the axis with $|a|$. In every case the end has finite arc length and $H\to\infty$.
\end{proof}

\begin{figure}[tp]
\centering
\includegraphics[width=0.72\textwidth]{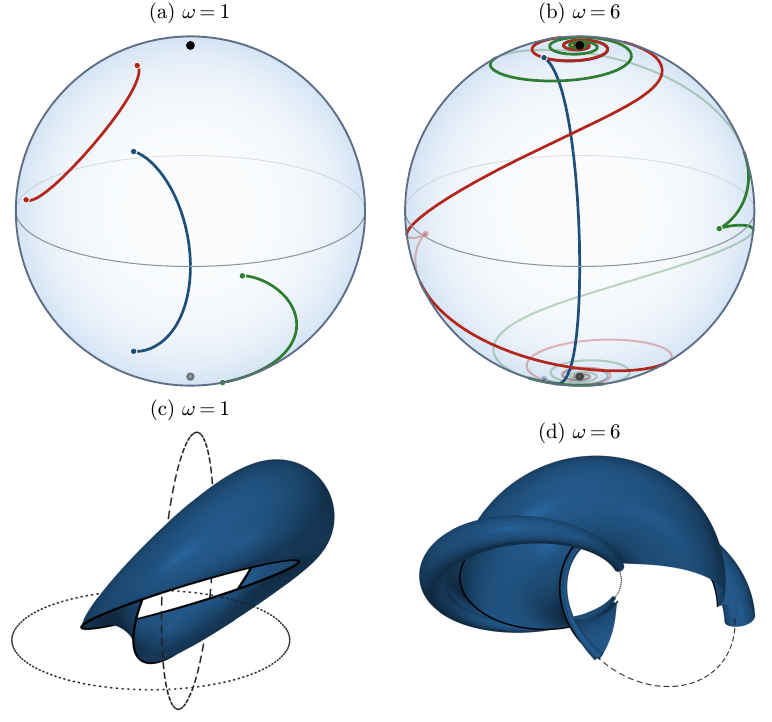}
\caption{Hopf solitons of $\Sp^3$ (Theorem~\ref{thm:hopf}). (a), (b) Rotators of the inverse curve shortening flow on $\Sp^2(\frac12)$; the black dots are the fixed points of the rotation, and the parts of the curves on the back of the sphere are drawn lighter. For $\omega=1<2$ each rotator ends, in both directions, at a point where $\kappa\to\infty$ (dots), the image of a singular circle; for $\omega=6\ge2$ the rotators reach the two fixed points of the rotation along logarithmic spirals. (c), (d) The Hopf-invariant surfaces over the blue curves, under a stereographic projection of $\Sp^3$. The dotted and dashed circles are the fibres over the two fixed points (the axis and the far axis), which form a Hopf link. In (c) the surface is an annulus bounded by two cuspidal edges, which are Hopf fibres; in (d) the surface is cut open, and it spirals towards both circles.}
\label{fig:hopf}
\end{figure}

\begin{proof}[Proof of Theorem~\ref{thm:hopf}]
In complex notation, $Y_x=ix$ and $|Y|\equiv1$. By \eqref{eq:quot}, $g_{\rm quot}=dr^2+\frac14\sin^2(2r)\,d\psi^2$, which is the round sphere $\Sp^2(\frac12)$, and $2H=\kappa$ is the geodesic curvature of the projected curve \cite{Pi}. The horizontal part of $X$ projects to $(p-q)\partial_\psi$, the rotation fixing the images of the two axes, and $\la X,N\ra=\la\bar X,\bar N\ra$. Hence \eqref{eq:soliton} is exactly the rotator equation $\kappa\la\frac12\bar X,\bar n\ra=1$, with angular speed $\omega=\frac12(p-q)$. Completeness is preserved, because the projection is a Riemannian submersion with compact fibres. The remaining assertions follow from Theorem~\ref{thm:helic}(ii) with $|a|=|b|=1$, since $\lambda\ge4$ means $\omega\ge2$.
\end{proof}

\section{Numerical illustrations}\label{sec:numerics}

None of the results above depends on the computations below. We include them as illustrations of the constants that appear in the statements. All integrations use RK45 (SciPy) with $\mathrm{rtol}=10^{-11}$ and $\mathrm{atol}=10^{-13}$, and the Radau method for the stiff problems of the horospherical ends. Caps are launched at small $r_0$ with $v_0=r_0/\lambda$, as in Lemma~\ref{lem:cap}.

\emph{Constants in dimension $n$.} The table below gives the Euclidean constant $c_0(n)$ of Theorem~\ref{thm:rplus}(i), computed as $r_+(\lambda)/\lambda$ for $\lambda=0.004$ (the values for $\Sp^n$ and $\Hy^n$ agree to $10^{-5}$), and compares $\lambda(\frac\pi2-r_+)$ at $\lambda=4000$ with $K_n$. The last two rows test the logarithmic term of Theorem~\ref{thm:rplus}(ii): the quantity
\[
Q_n(\lambda)=\lambda^2\Big(\lambda\big(\tfrac\pi2-r_+\big)-K_n\Big)+2(n-1)^2K_n\log\lambda
\]
is bounded by the theorem, and numerically it is constant to four digits between $\lambda=1000$ and $\lambda=4000$. For every $n$ in the table, bisection on the alternative ``entire versus singular sphere'' for the translators of $\Hy^n$ locates the critical speed $1$ to within $10^{-9}$, as predicted by Theorem~\ref{thm:transl}.
\begin{center}\small
\begin{tabular}{lcccccc}\toprule
$n$ & 3 & 4 & 5 & 6 & 8 & 10\\\midrule
$c_0(n)$ & 0.72215 & 0.68774 & 0.66526 & 0.64908 & 0.62687 & 0.61202\\
$\lambda(\frac\pi2-r_+)$, $\lambda=4000$ & 1.99999 & 2.35618 & 2.66664 & 2.94519 & 3.43600 & 3.86543\\
$K_n$ & 2 & 2.35619 & 2.66667 & 2.94524 & 3.43612 & 3.86563\\
$Q_n(1000)$ & 25.97 & 87.33 & 201.27 & 381.08 & 986.65 & 1993.8\\
$Q_n(4000)$ & 25.97 & 87.33 & 201.29 & 381.13 & 986.87 & 1994.5\\\bottomrule
\end{tabular}
\end{center}
For the hyperbolic type with $\lambda=4(n-1)$, inward solutions above the exceptional one reach the axis with $v\to0.965926=\cos\frac\pi{12}$ for $n=3,4,5$, the larger root of $v\sqrt{1-v^2}=\frac14$. Bisection between the two fates of Theorem~\ref{thm:hyp}(ii) locates the exceptional solution, for instance $v_e(0.3)=0.247328$ for $n=3$, and along it $v\to0.258819=\sin\frac\pi{12}$, the smaller root. In the borderline case $\lambda=2(n-1)$, the threshold solution $v_c$ of Theorem~\ref{thm:hyp}(iii) has $(v_c-\frac1{\sqrt2})/r^2=-0.5299$, $-0.8828$ and $-1.2353$ at $r=0.03$ for $n=3,4,5$, to be compared with $-\frac{2n-3}{4\sqrt2}=-0.5303$, $-0.8839$ and $-1.2374$. For the caps of $\Sp^4$ and $\Sp^5$ with $\lambda=3(n-1)$, $\sigma\to1.935660=\frac\pi2+\frac12\arcsin\frac23\pmod\pi$. Integrating backwards from $r=8$, the outer ends of both types in $\Hy^4$ and $\Hy^5$ with $\lambda=2$ have $z/z_*=1.0000$ and $H=1.0000$ at $r=5$, independently of the initial datum. Along the parabolic separatrix, $(n-1)\alpha/h=1.0000$ for $n=3,4,5$.

\emph{The case $n=3$.} The table below compares $r_+(\lambda)$ with the leading terms of Corollary~\ref{cor:rplus3}; the deficit $2-\lambda(\frac\pi2-r_+)$ is $(16\log\lambda-25.97)/\lambda^2$ up to higher order terms. As predicted by the dependence on $\varepsilon\lambda^2$, the second-order coefficients for $\Sp^3$ and $\Hy^3$ have the same size and opposite signs: $r_+/\lambda-c_0\approx-0.18\,\varepsilon\lambda^2$.
\begin{center}\small
\begin{tabular}{lcccc|lccccc}\toprule
$\lambda$ & 0.4 & 0.2 & 0.1 & 0.04 & $\lambda$ & 8 & 20 & 40 & 200 & 400\\\midrule
$r_+/\lambda$, $\Sp^3$ & 0.6957 & 0.7152 & 0.7204 & 0.7219 & $\lambda(\frac\pi2-r_+)$, $\Sp^3$ & 1.812 & 1.942 & 1.979 & 1.9985 & 1.9996\\
$r_+/\lambda$, $\Hy^3$ & 0.7533 & 0.7294 & 0.7240 & 0.7224 & & & & & &\\\bottomrule
\end{tabular}
\end{center}
In $\Hy^3$, $r_+(0.6)=0.4813$, $r_+(0.7)=0.5900$ and $r_+(0.9)=0.9135$, increasing to $\infty$ as $\lambda\uparrow1$. For the cap of $\Sp^3$ with $\lambda=8$, $\sigma\to1.832596=\sigma_*$, and the winding angle grows by $17.2=|\tan\sigma_*|\log100$ each time $x$ is divided by $100$. In the borderline case $\lambda=4$, the quantity $(\frac{3\pi}4-\sigma)|\log x|$ takes the values $0.29$, $0.36$, $0.39$ and $0.41$ at $x=10^{-3},10^{-6},10^{-9},10^{-12}$. This is the slow logarithmic convergence to $\frac12$ predicted by Lemma~\ref{lem:sn}. On the line $h=0.3$, bisection locates the parabolic separatrix at $\alpha_0=0.1505589$, and backward along it $\alpha/h$ stays equal to $0.50$; integrating forward from the centre manifold at $h=0.004$ gives the same value of $\alpha_0$ to seven digits.

\emph{Figures.} Figures~\ref{fig:transl}--\ref{fig:hopf}, placed next to the corresponding proofs, illustrate the case $n=3$. They show solutions of the reduced equations for each class, together with the corresponding surfaces. The profiles were obtained by integrating the system \eqref{eq:systemt} in the parameter $t$ with $ds=\cos\sigma\,dt$, which is regular across the singular orbits, so that each curve is drawn as a front; dots mark the cuspidal points. The surfaces are drawn through the stereographic projection $\pi(x)=(x_1,x_2,x_3)/(1+x_0)$, which maps $\Hy^3$ onto the Poincar\'e ball and $\Sp^3\setminus\{-e_0\}$ onto $\R^3$. The ideal boundary of $\Hy^3$ is drawn as a transparent sphere, and black curves on the surfaces are cuspidal edges.

\subsection*{Acknowledgements}
The authors thank the Instituto de Matem\'atica e Estat\'istica of the Universidade Federal de Goi\'as, where this work was carried out.

\subsection*{Declarations}
\textbf{Funding.} This research received no specific funding.
\textbf{Conflict of interest.} The authors declare that they have no conflict of interest.
\textbf{Use of AI tools.} The authors used the AI assistant Claude (Anthropic) to assist with language editing, typesetting, figure generation, and numerical verification. All AI-generated content and results were carefully reviewed and verified by the authors, who take full responsibility for the final manuscript.
\textbf{Data availability.} No datasets were generated or analysed in this work. The scripts used for the numerical illustrations and the figures of Section~\ref{sec:numerics} are available from the corresponding author upon request.

\end{document}